\documentclass[11pt,twoside]{article}

\usepackage{graphicx}
\usepackage{multirow}
\usepackage{fancyhdr}
\usepackage{array}
\usepackage{amsfonts}
\usepackage{amsmath}
\usepackage{amssymb}
\usepackage{amsthm}

\usepackage{url}
\usepackage{algorithm}
\usepackage[noend]{algpseudocode}
\makeatletter
\def\BState{\State\hskip-\ALG@thistlm}
\makeatother

\usepackage{subfigure}
\usepackage{color}

\usepackage[textwidth=16cm,textheight=22cm,top=3cm,bottom=3cm]{geometry}

\newtheorem{theorem}{Theorem}
\newtheorem{remark}{Remark}
\newtheorem{corollary}{Corollary}

\graphicspath{{./figures/}}

\newcommand{\drop}[1]{}
\newcommand{\no}{\noindent}
\newcommand{\pt}{{\partial_t}}
\newcommand{\p}{{\partial}}

\newcommand{\vfi}{\varphi}
\newcommand{\R}{\mathbb{R}}
\newcommand{\N}{\mathbb{N}}
\newcommand{\C}{\mathbb{C}}

\newcommand{\fer}[1]{(\ref{#1})}
\newcommand{\qtext}[1]{\quad\text{#1}}
\newcommand{\abs}[1]{| #1 |}
\newcommand{\nor}[1]{\| #1 \|}
\newcommand{\grad}{\nabla}
\def\O{\Omega}

\newcommand{\bu}{\mathbf{u}}
\newcommand{\bv}{\mathbf{v}}

\newcommand{\bx}{\mathbf{x}}
\newcommand{\by}{\mathbf{y}}
\newcommand{\bz}{\mathbf{z}}

\newcommand{\wto}{\rightharpoonup}

\DeclareMathOperator{\Div}{div}
\DeclareMathOperator{\tr}{tr}

\DeclareMathOperator{\diag}{diag}

\DeclareMathOperator{\SNR}{SNR}

\DeclareMathOperator{\Ya}{Y}
\DeclareMathOperator{\NLM}{NLM}

\title{ On a cross-diffusion system arising in image denosing 
\thanks{Supported by Spanish MCI Project MTM2017-87162-P.}}

\author{Gonzalo Galiano  \thanks{Dpt. of Mathematics, Universidad de Oviedo,
 c/ Calvo Sotelo, 33007-Oviedo, Spain ({\tt galiano@uniovi.es, julian@uniovi.es})}
    \and Juli\'an Velasco\footnotemark[2] }

\date{}

\begin{document}

\thispagestyle{plain}

\maketitle

\begin{abstract}
We study a generalization of a cross-diffusion problem deduced from a nonlinear complex-variable 
diffusion model for signal and image denoising. 

We prove the existence of weak solutions of the time-independent problem with fidelity terms under mild conditions on the data problem. Then, we show that this translates on the well-posedness of a   
quasi-steady state approximation of the evolution problem, and also prove the existence 
of weak solutions of the latter under more restrictive hypothesis.

We finally perform some numerical simulations for image denoising, comparing the performance of the cross-diffusion model to other well-known demoising methods.

\no\emph{Keywords: }cross-diffusion, denoising, Perona-Malik, existence of solutions, simulations.

\end{abstract}

\section{Introduction}

A fundamental problem in image analysis is the reconstruction of an original image from an observation,
usually to perform subsequent operations on it, such as segmentation or inpainting.
In this reconstruction, image denoising is the main task, involving the elimination or reduction 
of random phenomena -noise- which has been introduced in the image acquisition process due to a number of
factors like poor lighting, magnetic perturbances, etc.

Additive white Gaussian noise is surely the most used noise model. In it, we assume that the observation 
$u_0$ is the addition of a \emph{ground truth} image, $u$, with a zero mean Gaussian random variable, $n$. That is, $u_0 = u + n$. The problem is then how to recover $u$ from $u_0$.

Along the years, many different techniques have been developed to treat with this
important technological problem. Just to mention a few well known approaches, we have from variational models such 
as the ROF model \cite{Rudin1992}, to models based on neighborhoods such as the Bilateral filters \cite{Yaroslavsky1985, Smith1997, Tomasi1998}, 
based on nonlocal information such as the Nonlocal Means \cite{Buades2005}, or to methods based on partial differential
equations, such as anisotropic diffusion models \cite{Alvarez1992}, or the Perona-Malik equation \cite{Perona1990}.

The Perona-Malik equation has received much attention due to the confluence of the very good denoising 
performance on the discrete level, the mathematical simplicity of its formulation, and the enormous 
difficulties to establish the well-posedness of the mathematical problem \cite{Guidotti2014}. It reads as follows:
\begin{align}
\label{eq:pm}
 \p_t u = \Div( g\grad u ). 
\end{align}
The function $g:\R_+\to\R_+$ is known as the \emph{edge detector}, which prevents or induces diffusion according 
to some criteria. In the most usual model, $g \equiv g(\abs{\grad u})$, with $g$ satisfying $g(0) >0$ and $g(s)\to 0$ as $s\to\infty$. Common examples of the edge detector are 
\begin{align}
 \label{def:g}
 g(s) = \frac{1}{1+(s/\lambda)^2},\qtext{and}\quad 
 g(s) = \exp\big( -\frac{s^2}{\lambda^2}\big),
\end{align}
for some scaling parameter $\lambda >0 $. 
Thus, for large gradients (edges in the image) the diffusion
is small and so the edges are preserved. On the contrary, for small gradients (almost flat regions) 
diffusion is allowed, contributing to the elimination of small perturbances, i.e. denoising the image.
Similar results are obtained when the norm of the gradient is replaced by the absolute value of 
the Laplacian, that is, when considering the edge detector as $g(\abs{\Delta u})$.

As aforementioned, many efforts have been devoted to prove the well-posedness of the Perona-Malik 
equation, although many questions remain open. Other efforts have been directed to state well-posed
Perona-Malik type problems. In this line, Gilboa et al. \cite{Gilboa2001, Gilboa2004} 
introduced the following complex-variable problem: Find $I:[0,T]\times\O \to \C$ such that 
\begin{align*}
 &\p_t I = c\Div\big( g \grad I\big)&& \text{in } Q_T=(0,T)\times\O,&\\
 & g \grad I \cdot \nu = 0 && \text{on } \Gamma_T=(0,T)\times\p\O,&\\
 &I(0,\cdot) = I_0 && \text{in } \O,&
\end{align*}
where $\O\subset \R^m$ is a bounded smooth domain, $\p\O$ its boundary, $\nu$ the  exterior unitary normal vector, 
and $c = \exp(i\theta)$, for some $\theta\in(0,\pi/2)$. 
Introducing the decomposition of $I$ as the sum of its real and imaginary parts, $I= u_1+iu_2$, and expanding
the above equations, they arrived to the following problem: Find $u_1,u_2:[0,T]\times\O \to \R$ such that
\begin{align}
 &\p_t u_1 - \Div\big( g (a\grad u_1 - b \grad u_2) \big)=0&& \text{in } Q_T,& \label{cd:1}\\
 &\p_t u_2 - \Div\big( g (b\grad u_1 + a \grad u_2) \big)=0&& \text{in } Q_T,&\\
 & g \grad u_1 \cdot \nu = g \grad u_2 \cdot \nu = 0 && \text{on } \Gamma_T,&\\
 &u_1(0,\cdot) = {\rm Re}(I_0), ~~ u_2(0,\cdot) = {\rm Im}(I_0), && \text{in } \O,& \label{cd:2}
\end{align}
with $a = \cos\theta$, $b=\sin\theta$. In this denoising model, $I_0$ is identified with the initial 
noisy image, leading to ${\rm Re}(I_0) = u_0$, and ${\rm Im}(I_0) = 0$. Thus, the $u_1$ equation 
establishes the rules of evolution of the noisy image, so $u_1$ should be regarded as an image 
itself. The second equation may be rewritten as 
\begin{align}
\label{eq:key}
 \p_t u_2  - a \Div\big( g  \grad u_2 \big) = b \Div\big( g \grad u_1  \big),
\end{align}
showing that $u_2$ is related to some regularization involving the second order derivatives of $u_1$. This motivated Gilboa et al. to, instead of considering an edge detector of the form $g(\abs{\Delta u_1})$, 
use the detector $g \equiv g(u_2)$.

The Perona-Malik equation \fer{eq:pm} is the pure filtering problem for denoising the original image. When deduced from variational principles, this kind of equations arise as the gradient descent corresponding to the Euler-Lagrange equation satisfied by the minima of certain functional. In this variational framework, evolution equations such as 
\begin{align*}
 \p_t u -\Div(g\grad u) = \beta (u-u_0),
\end{align*}
may be used to approximate the \emph{true} filtered image, solution of
\begin{align*}
 -\Div(g\grad u) = \beta (u-u_0).
\end{align*}
Here, the right hand side is the \emph{fidelity term}, introduced to avoid an excessive smoothing in the denoising procedure, and $\beta >0$ is a constant balancing regularization  and fidelity. 

The corresponding extension of Gilboa's et al. model \fer{cd:1}-\fer{cd:2} to its steady state version is: Find 
$u_1,u_2:\O \to \R$ such that
\begin{align}
 &- \Div\big( g (a\grad u_1 - b \grad u_2) \big)=\beta_1(u_{10}-u_1)&& \text{in } \O,& \label{cds:1}\\
 & - \Div\big( g (b\grad u_1 + a \grad u_2) \big)=\beta_2(u_{20}-u_2)&& \text{in } \O,&\\
 & g \grad u_1 \cdot \nu = g \grad u_2 \cdot \nu = 0 && \text{on } \p\O,& \label{cds:2}
\end{align}
being $u_{10}$ the original noisy image, and $u_{20}=0$.

The proof of existence of solutions of problem \fer{cd:1}-\fer{cd:2} (with or without fidelity terms) and of problem \fer{cds:1}-\fer{cds:2} may be achieved under several sets of restrictions on the data. We mention here three of such cases, making reference only to the evolution problem. The changes which must be introduced to cover the steady-state problem are straightforward.

\noindent\emph{Case 1. } When $g$ is a positive constant, so that the problem becomes linear. Then, if the resulting diffusion matrix, which may be more general than that in \fer{cd:1}-\fer{cd:2}, is definite positive the existence of solutions is a classical result.  

\noindent\emph{Case 2. }
Under the following restriction on the edge detector: there exists a positive constant $\alpha$ such that $\alpha \leq g(s) \leq \alpha^{-1}$, for all $s\in \R_+$. In this case, the following key a priori estimate holds for the standard energies 
\begin{align}
\label{energy}
\sum_{i=1}^2 \Big(\int_\O u_i^2 + \alpha \int_{Q_T} \abs{\grad u_i}^2\Big) \leq \int_\O u_{i0}^2,
\end{align}
providing the compactness needed to construct the solution from approximation arguments.

\noindent\emph{Case 3. } When $b=0$ the cross-diffusion is eliminated from the problem. Then,  $L^\infty$ estimates for $u_1,~u_2$ are obtained by standard arguments, making estimate \fer{energy} valid.
In this line, one can ask whether a change of unknowns in \fer{cd:1}-\fer{cd:2} may lead to a diagonal diffusion system.  Assuming the following general form for the diffusion terms 
\begin{align*}
 - \Div\big( g(u_1,u_2) (a_{i1}\grad u_1 +a_{i2} \grad u_2) \big)
\end{align*}
one finds that this change is possible, and therefore existence of solutions is granted, if one of the following conditions is satisfied:
\begin{enumerate}
 \item $a_{12}a_{21} >0$, or
 \item $a_{12}a_{21} = 0$ and $a_{11}\neq a_{22}$,
\end{enumerate}
Unfortunately, none of them is satisfied by problem \fer{cd:1}-\fer{cd:2}, although the second is on the base of the proof of existence of solutions given in \cite{Lorenz2006} for a 
simplified triangular version of the problem, already proposed by Gilboa et al. \cite{Gilboa2001}.

All these special cases have some disadvantages from the image denoising point of view. The linear filtering, as it is well known, produces an over-smoothing which results in the  edges of the objects within the image being blurred. In the diagonal or triangular matrix diffusion case, the effects of cross-diffusion are neglected or very mild, and the interpretation of the second component, $u_2$, as a regularized version of $\Delta u_1$, see \fer{eq:key}, does not apply. Finally, considering an edge detector which is a priori bounded away from zero introduces a new parameter, $\alpha$, which must be estimated for a practical implementation of the filter. 

In this article we study the model problems \fer{cd:1}-\fer{cd:2} and \fer{cds:1}-\fer{cds:2}
without any additional assumptions on the edge detector, $g$, and for a more general form of the diffusion coefficients. Moreover, the possibility of including a fidelity term in the evolution problem  \fer{cd:1}-\fer{cd:2} is also treated.

The main result is the proof of existence of weak solutions of the steady state problem, see Theorem~\ref{th:steady}. This result is instrumental for establishing the well-posedness of a quasi-stationary steady state (QSS) approximation of the (generalized) evolution problem \fer{cd:1}-\fer{cd:2}, see Corollary~\ref{cor:qssa}. The QSS can be interpreted as a discrete-time gradient descent approximation to the solution of \fer{cd:1}-\fer{cd:2}, which we use in Section~\ref{sec:numerics} to illustrate the denoising features of the evolution cross-diffusion model.  

However, our $L^\infty$ estimates for the QSS approximation depend, in general, on the time discretization step, $\tau$, preventing the passing to the limit $\tau\to 0$ to get  a solution of the time continuous  problem. In Corollary~\ref{cor:diagonal}, we summarize the special cases in which an uniform $L^\infty$ estimate may be found, and thus the existence of solutions of the evolution problem may be proven. 

Let us finally mention that cross-diffusion problems like \fer{cd:1}-\fer{cd:2} appear in many applications in physics, chemistry or biology for describing the time evolution of densities of multi-component systems. 
A fundamental theory for the study of strongly coupled
systems was developed by Amann \cite{Amann1989} which established, under suitable conditions, local existence of solutions, which become global if their $L^\infty$ 
and H\"older norms may be controlled. 

Since generally no maximum
principle holds for parabolic systems, the proof of $L^\infty$ bounds is a challenging problem \cite{Jungel2016}, and ad hoc techniques must be employed 
to deduce the existence of solutions of particular problems. A useful methodology avoiding the $L^\infty$ bounds requirement of Amman's results was introduced in \cite{Galiano2003} and later generalized in a series of papers, see  \cite{Chen2004,Chen2006,Galiano2012,Desvillettes2014,Galiano2014,Jungel2016} and the references therein. However, this technique relies on the introduction of a Lyapunov functional  needing of a particular cross-diffusion structure which is not satisfied by system \fer{cd:1}-\fer{cd:2}. 
Thus, the proof of existence of solutions of problem \fer{cd:1}-\fer{cd:2} with the edge detector $g$ of the type \fer{def:g} remains open.

The nonlinear instability of these type of strongly coupled systems, as induced only by cross diffusion terms, has been also an active area of research, see  \cite{Gambino2013, Ruiz2013, Cai2016, Gambino2018}  and the references therein. The investigation of this property for systems of the type \fer{cd:1}-\fer{cd:2} will be subject of future research.

\section{Main results}

We consider the following problems.

\smallskip

\no\emph{Evolution problem: } Given a fixed $T>0$ and a bounded domain $\O\subset\R^m$,  find $\bu=(u_1,u_2)$, with 
$u_i:(0,T)\times\O\to\R$, such that, for $i=1,2$,
\begin{align}
& \pt u_i-\Div J_i(\bu)=f_i(\bu) && \qtext{in }Q_T, 
  \label{eq:pde}\\
& J_i(\bu)\cdot \nu =0 && \qtext{on }\Gamma_T,
\label{eq:bc}\\
& u_{i}(\cdot,0)=u_{i0} && \qtext{in }\O,	\label{eq:id} 
\end{align}
with flow and fidelity functions given by
\begin{align}
 & J_i(\bu) = g(\bu)\big(a_{i1} \grad u_1 + a_{i2} \grad u_2\big),  \label{def:flow}\\
& f_i(\bu) = \beta_i(u_{i0}-u_i) ,
\label{def:reaction}
\end{align}
for $\bu_0 \in L^\infty(\O)^2$, and for some constants $\beta_i \geq 0$.
 
\smallskip 
 
\no\emph{Steady state problem: } The steady state corresponding to problem \fer{eq:pde}-\fer{eq:id} is found by dropping the time derivative in \fer{eq:pde}. However, we shall study the following generalization: find $\bu=(u_1,u_2)$, with 
$u_i:\O\to\R$, such that, for $i=1,2$,
\begin{align}
& \gamma_iu_i -\Div J_i(\bu)=G_i  && \qtext{in }\O, 
  \label{eqs:pde}\\
& J_i(\bu)\cdot \nu =0 && \qtext{on }\p\O,
\label{eqs:bc} 
\end{align}
for some functions $G_i\in L^\infty(\O)$, constants $\gamma_i > 0$, and $J_i$ given by \fer{def:flow}.

\smallskip 

\no\emph{Assumptions on the data. } We make the following hypothesis on the data, which we shall refer to as \textbf{(H)}:
\begin{enumerate}
\item $\O\subset\R^m$ is a bounded domain with Lipschitz continuous boundary, $\partial \O$.

\item The constant matrix $A= (a_{ij})$ satisfies:
\begin{itemize}
 \item There exists a constant $a_0 >0$ such that 
\begin{equation*}
\xi^T A~\xi \geq a_0 \abs{\xi}^2 \qtext{for all }\xi\in\R^m. 
\end{equation*}

\item The elements of $A$ satisfy $a_{ii}> \abs{a_{ji}}$ for $i,j=1,2$, $i\neq j$. 
\end{itemize}

\item  The \emph{edge detector} $g:\R^2\to\R$ is continuous, with $g(K)\subset(0,\infty)$ for all compact $K\subset\R^2$.

\end{enumerate}

\begin{theorem}\label{th:steady}
Assume (H), and let $\gamma_i >0$, and  $G_i\in L^\infty(\O)$, for $i=1,2$. Then, there exists a weak solution  $\bu \in H^1(\O)^2\cap L^\infty(\O)^2$,
of the stationary problem \fer{eqs:pde}-\fer{eqs:bc}.
\end{theorem}
The next result establishes the existence of solutions of a time-discrete version of 
problem \fer{eq:pde}-\fer{eq:id} known as the \emph{quasi-steady state approximation} 
problem. This is the time discretization scheme we use in our numerical experiments of  Section~\ref{sec:numerics}.

\begin{corollary}\label{cor:qssa}
Let $\tau>0$, assume (H), and set $ \bu^0=\bu_{0}\in L^\infty(\O)^2$. Then, for all $n\in\N$, there exists a weak solution $\bu^{n+1} = (u_1^{n+1},u_2^{n+1})$, with $u_i^{n+1} \in H^1(\O)\cap L^\infty(\O)$, of the problem,  for $i=1,2$, 
\begin{align}
& \frac{1}{\tau}(u_i^{n+1}- u_i^n )- \Div J_i(\bu^{n+1})=  f_i(\bu^{n+1}) && \qtext{in }\O, 
  \label{eq:pde:ss}\\
& J_i(\bu^{n+1})\cdot \nu =0 && \qtext{on }\partial\O.\label{eq:bc:ss}
\end{align}
\end{corollary}

Finally, from the construction in Theorem~\ref{th:steady} , we deduce  the existence of solutions of the evolution problem 
\fer{eq:pde}-\fer{eq:id} under some parameter restrictions.

\begin{corollary}\label{cor:diagonal}
Assume (H) and suppose that one of the following conditions is satisfied:
\begin{align*}
 (i)~ a_{12}=a_{21}=0,\quad (ii)~ a_{12}a_{21}>0, \text{ or } (iii)~ a_{12}a_{21}=0 \text{ and }a_{11}\neq a_{22}.
\end{align*}
Then, there exists a weak 
solution  $\bu = (u_1,u_2)$ of the evolution problem  \fer{eq:pde}-\fer{eq:id}, with, for 
$i=1,2$, 
\begin{align*}
u_i \in L^2(0,T;H^1(\O))\cap H^1(0,T; (H^1(\O))') \cap L^\infty(Q_T).
\end{align*}
\end{corollary}

\section{Numerical simulations}\label{sec:numerics}

In this section we numerically demonstrate the image denoising ability of the evolution cross-diffusion model \fer{eq:pde}-\fer{def:reaction} in its particular form \fer{cd:1}-\fer{cd:2} considered in Gilboa et al. \cite{Gilboa2001}. For comparison, we use two sets of methods. The first is based on the Perona-Malik (PM) equation which is, as already mentioned, closely related to the cross-diffusion system. With this comparison, we intend to check whether the cross-diffusion system is a clear improvement of the Perona-Malik equation or not. 

The second set is based on one step integral filters, among which we have chosen the Bilateral filter (BF) \cite{Yaroslavsky1985, Smith1997, Tomasi1998} and the Nonlocal Means filter (NLM) \cite{Buades2005}. The BF is well-known due to its simplicity and good denoising properties when compared to execution times. The NLM is a clear improvement of the BF regarding denoising quality, and it is one of the state of the art reference method to any denoising algorithm. Comparing to it we also indirectly compare to many other methods which take the NLM as  reference \cite{Buades2005}.

\begin{figure}[t]
\centering
 {\includegraphics[width=0.30\textwidth,height=0.21\textheight]{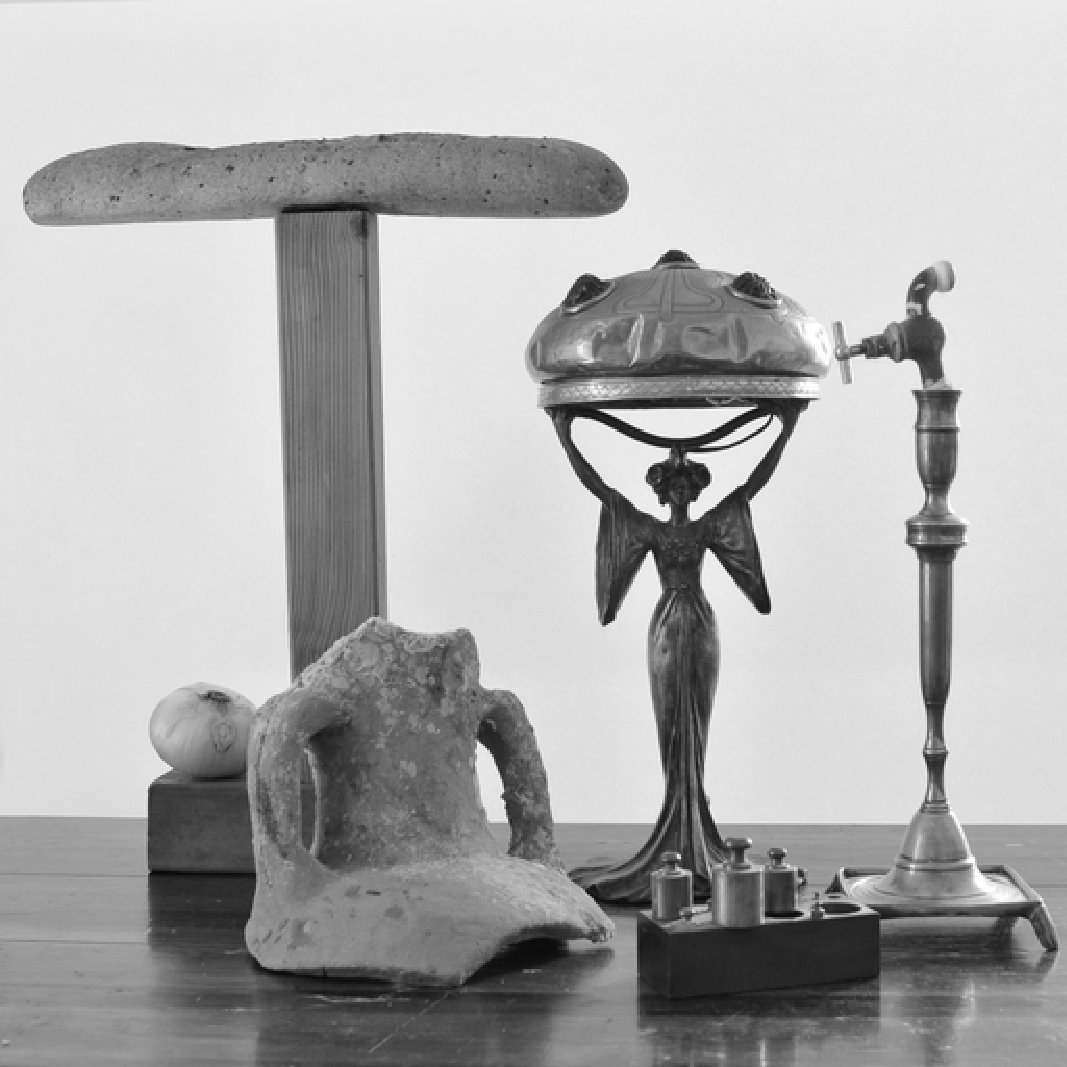}}
  {\includegraphics[width=0.30\textwidth,height=0.21\textheight]{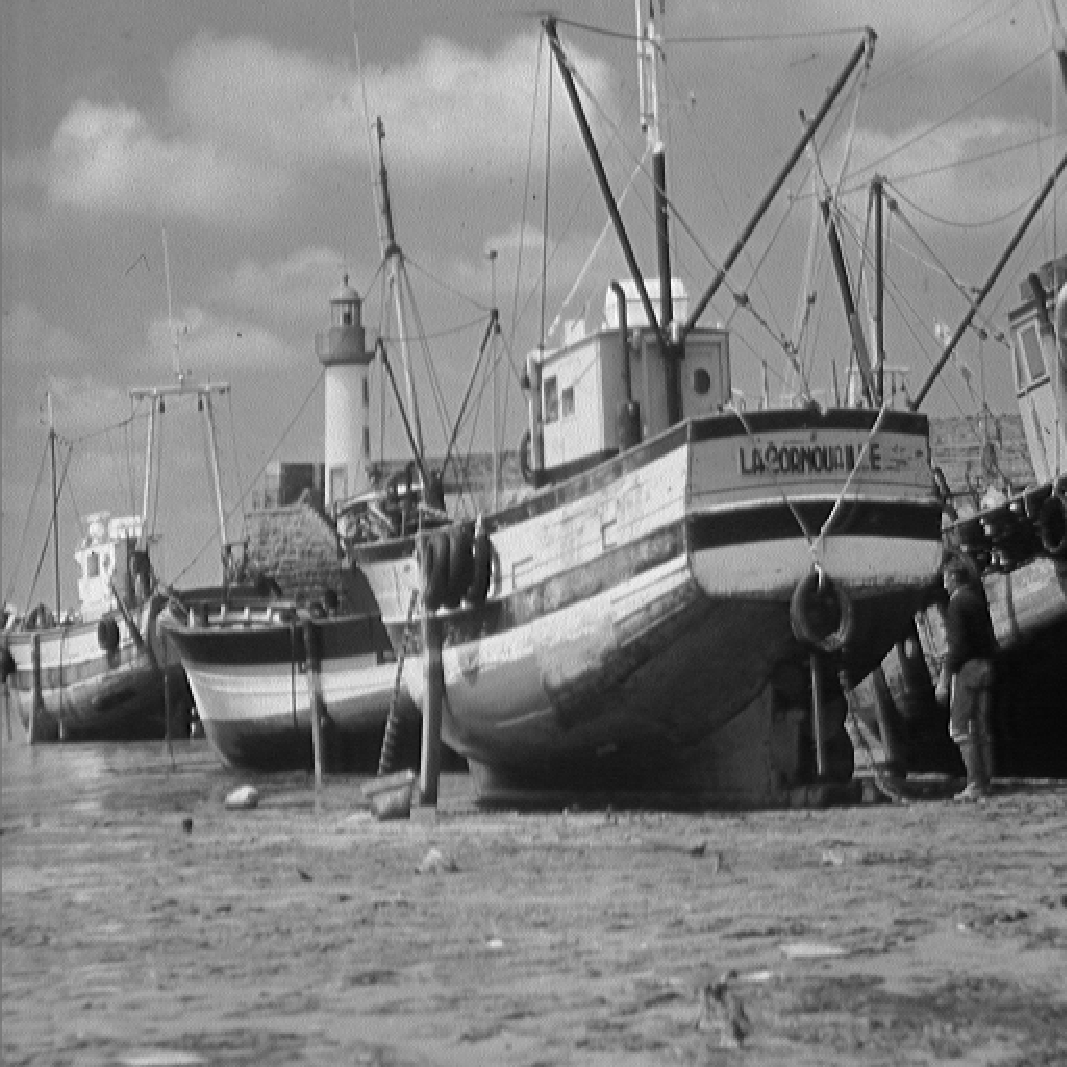}}
 {\includegraphics[width=0.30\textwidth,height=0.21\textheight]{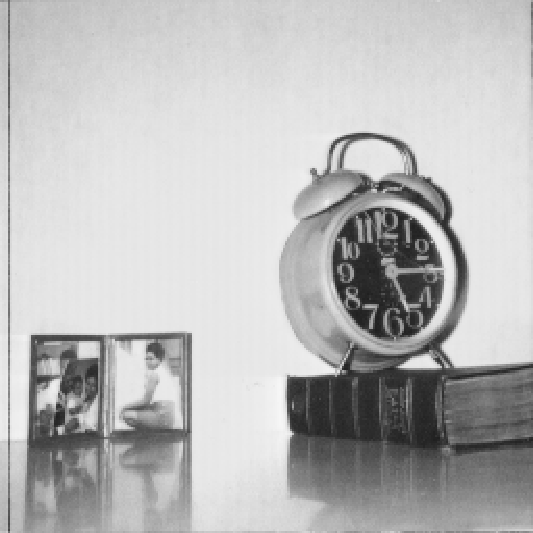}}\\
 {\includegraphics[width=0.30\textwidth,height=0.21\textheight]{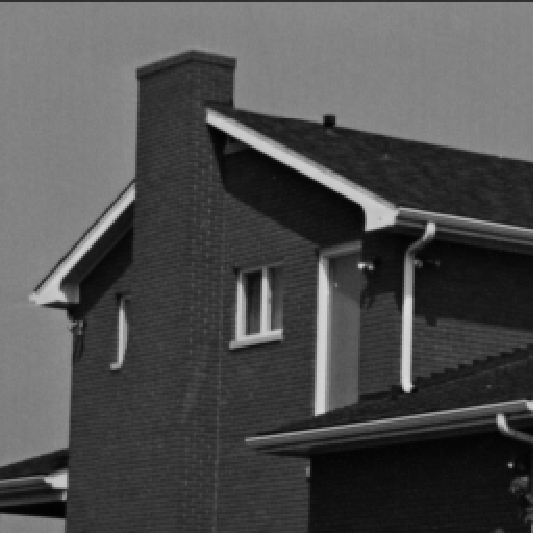}}
 {\includegraphics[width=0.30\textwidth,height=0.21\textheight]{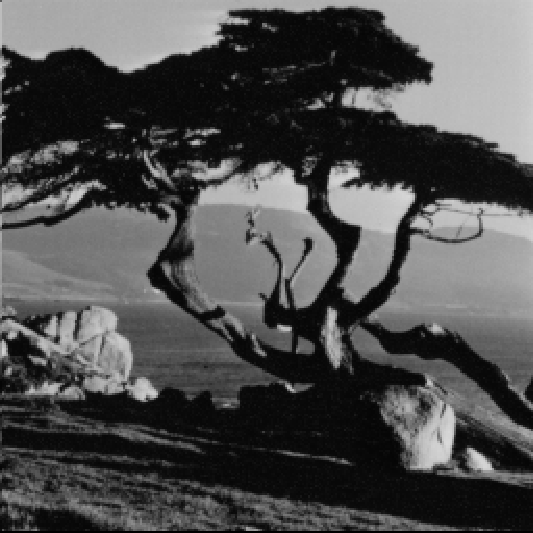}}
 {\includegraphics[width=0.30\textwidth,height=0.21\textheight]{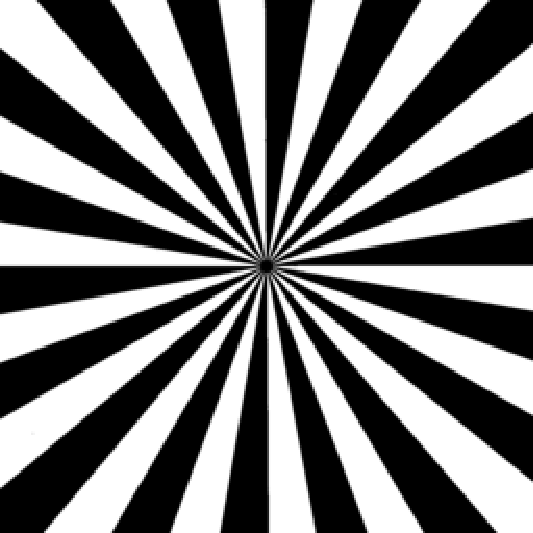}} \\
 {\includegraphics[width=0.18\textwidth,height=0.12\textheight]{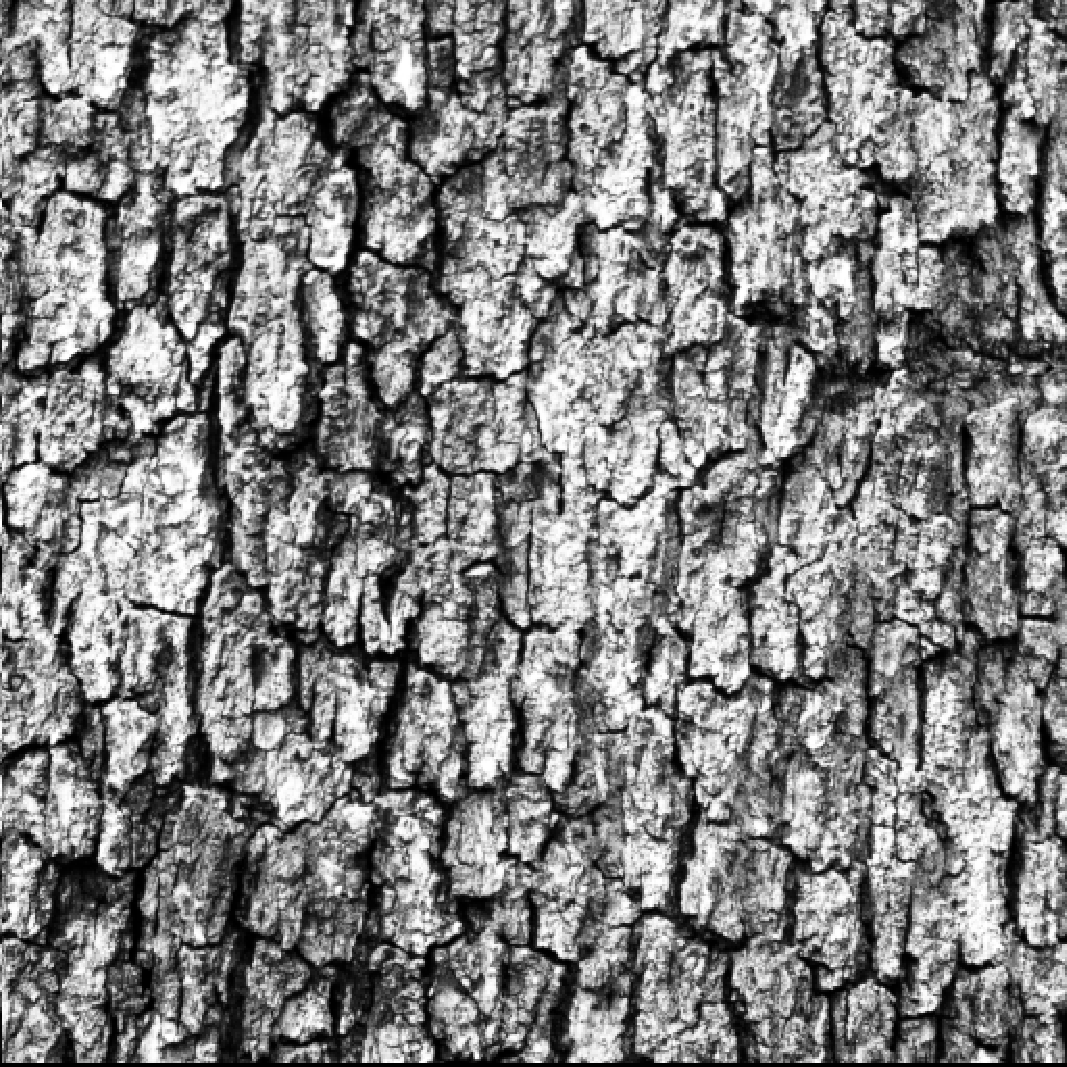}} 
 {\includegraphics[width=0.18\textwidth,height=0.12\textheight]{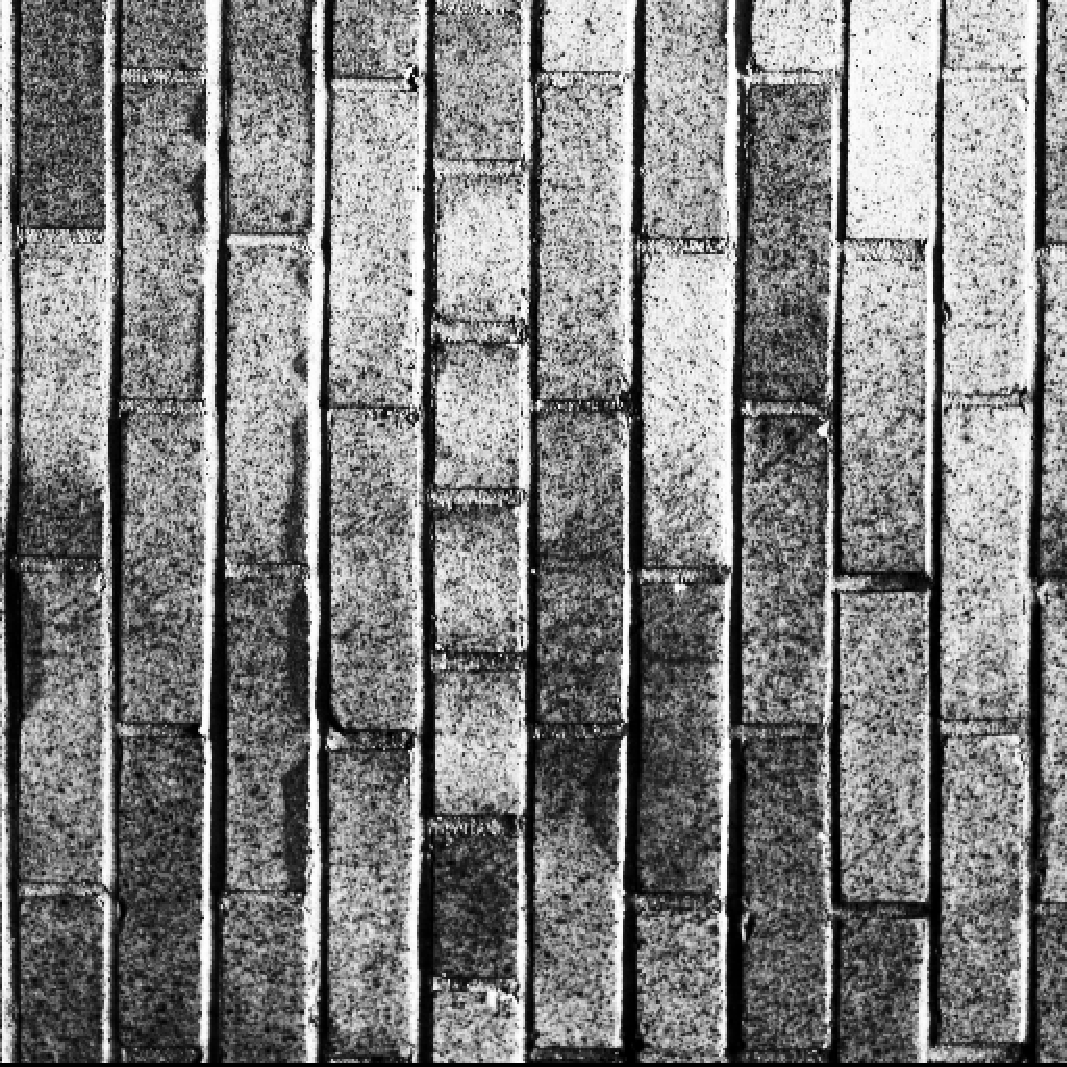}}  
 {\includegraphics[width=0.18\textwidth,height=0.12\textheight]{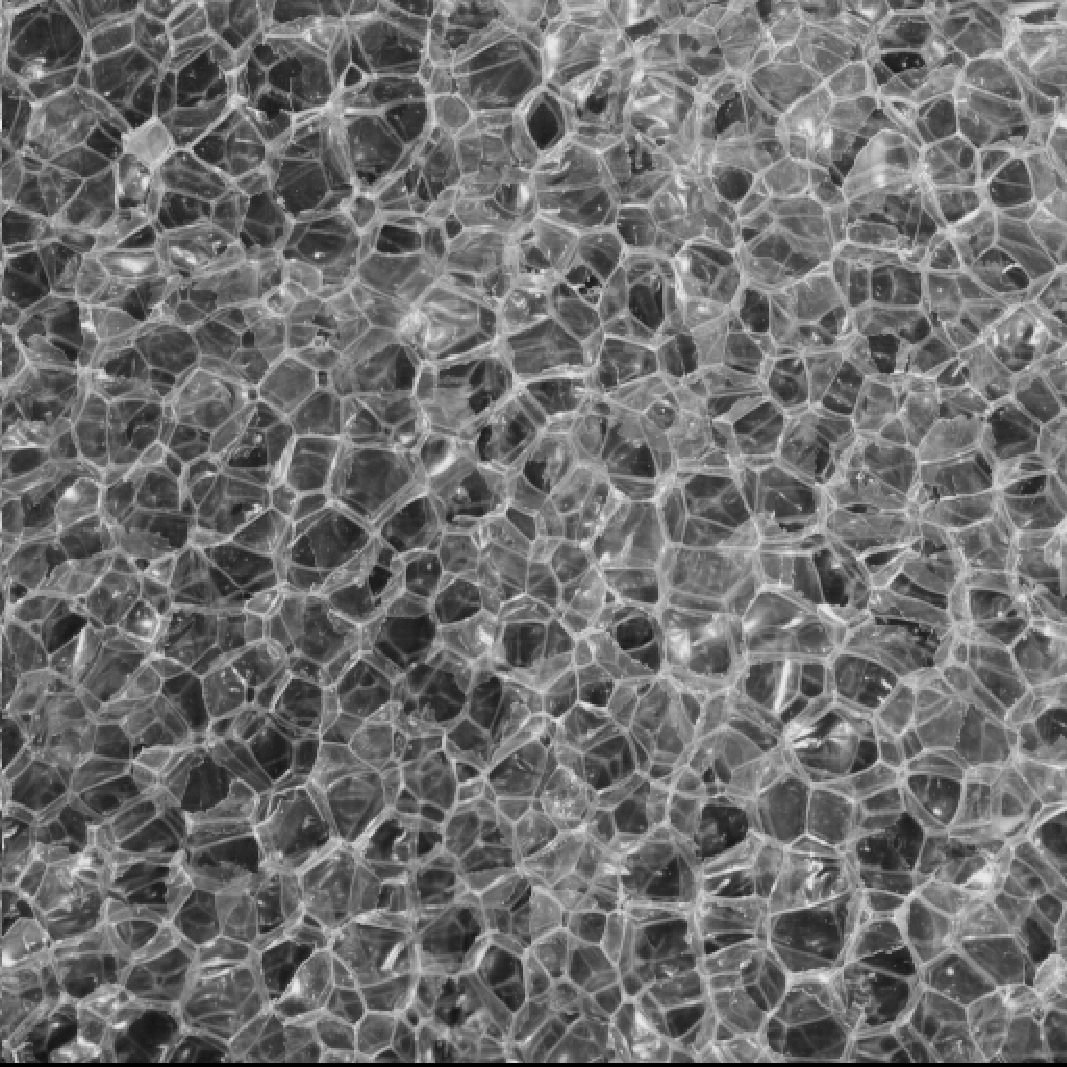}}
 {\includegraphics[width=0.18\textwidth,height=0.12\textheight]{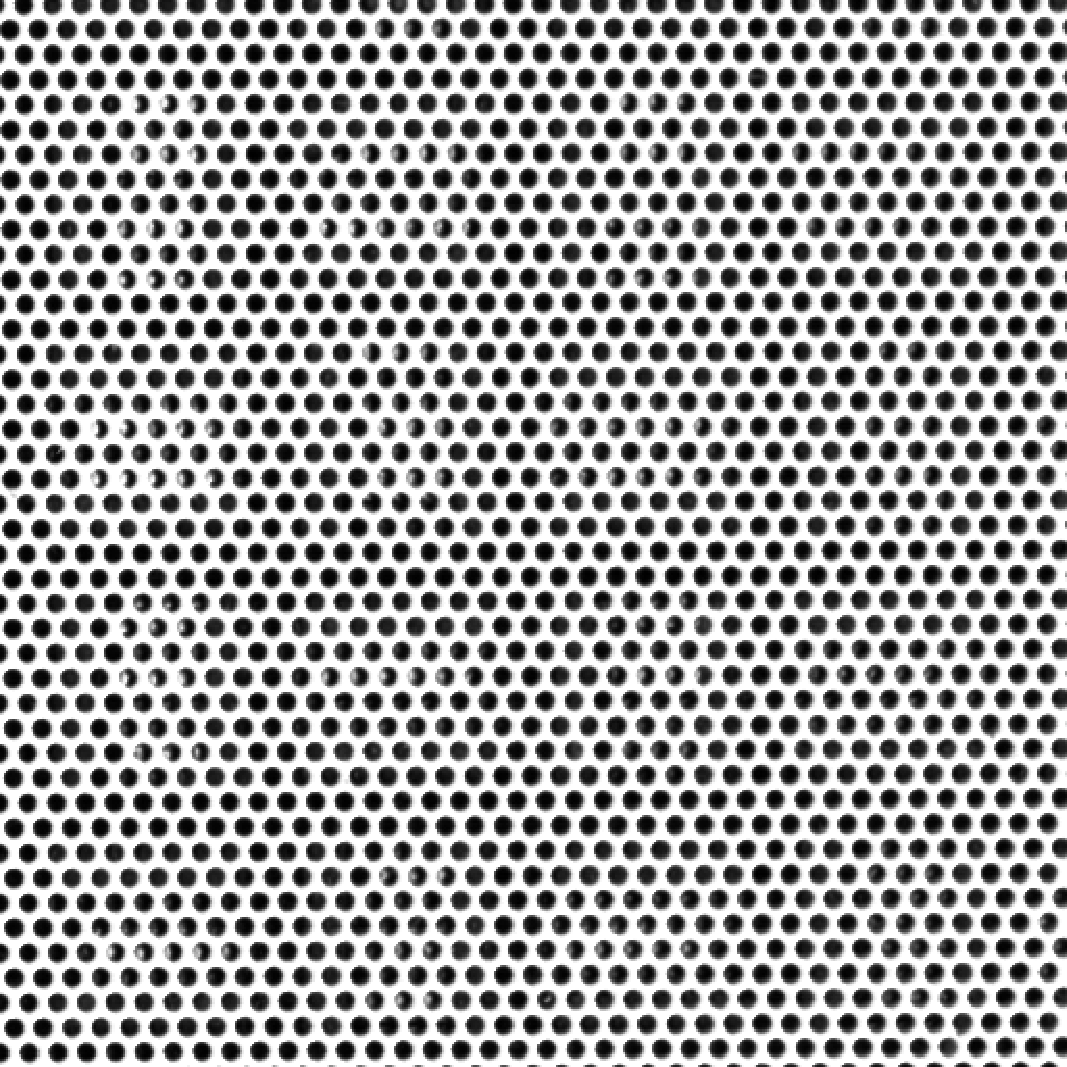}}
 {\includegraphics[width=0.18\textwidth,height=0.12\textheight]{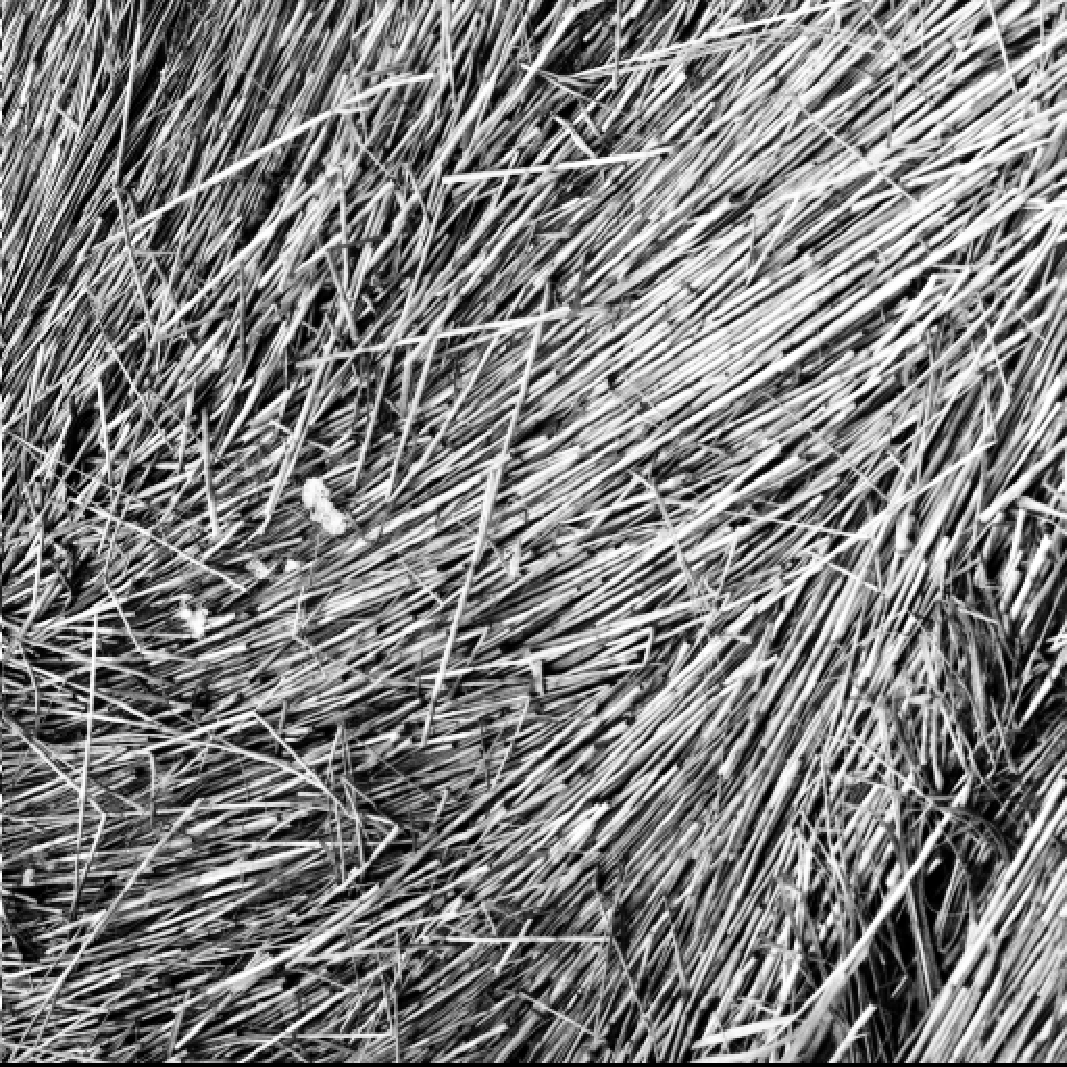}}
 \caption{Set of test images. First and second rows are natural images: \emph{added value, boat, clock, house, tree}, and \emph{test}. The latter is not \emph{natural} but it actually behaves as such in the denoising procedure. Third row are texture images: \emph{bark, bricks, bubbles,  holes}, and \emph{straw}. All of them are of size $512\times512$ but clock and test, which are $256\times256$.  } 
\label{fig1}
\end{figure}

The comparison of the methods is done on two classes of intensity images: natural images and texture images, most of them obtained from the SIPI\footnote{Signal and Image Processing Institute, University of Southern California.} data set, see Figure~\ref{fig1}. We add to the original image, $u_{\rm clean}$, a Gaussian white noise to obtain the noisy image, $u_{\rm noisy}$,  with fixed
signal to noise ratio $\SNR=10$, where  
\begin{equation*}
 \SNR = \frac{\sigma (u_{\rm clean})}{\sigma(u_{\rm clean}- u_{\rm noisy})},
\end{equation*}
being  $\sigma$ the standard deviation. After filtering $u_{\rm noisy}$ with each method, we obtain a denoised image, $u_{\rm denoised}$, that we compare to $u_{\rm clean}$
using three quality measures \cite{Wang2004}:
\begin{itemize}
 \item The pick signal to noise ratio (PSNR), given by 
\begin{align*}
 {\rm PSNR}(u_1,u_2) = 20 \log_{10} \Big(\frac{255}{\nor{u_1-u_2}_2}\Big).
\end{align*}

\item The normalized cross correlation (NCC), given by 
\begin{align*}
 {\rm NCC}(u_1,u_2) = \frac{\langle u_1, u_2 \rangle}{\nor{u_1}_2\nor{u_2}_2} \in [0,1].
\end{align*}

\item The structural similarity (SSIM), given by
\begin{align*}
 {\rm SSIM}(u_1,u_2) = \frac{(2\mu_1\mu_2+c_1)(2\sigma_{1,2}+c_2)}{(\mu_1^2+\mu_2^2+c_1)(\sigma_1^2+\sigma_2^2+c_2)} \in [0,1],
\end{align*}
where $\mu_i,~\sigma_i,~\sigma_{i,j}$ stand for the mean, the standard deviation and the covariance of the corresponding images, and $c_i$ are some positive small constants included to avoid instability when the denominators are close
to zero.
 
\end{itemize}

Since different methods have different sets of structural parameters which must be fixed for actual implementation, we have experimentally optimized them for each image in a reasonable range of values. The optimization has been performed with respect to the PSNR between the original and the denoised images.

\subsection{Discretization}

We use the well-posed QSS aproximation \fer{eq:pde:ss}-\fer{eq:bc:ss} of the cross-diffusion problem \fer{eq:pde}-\fer{eq:id} to discretize in time. A similar approach is followed for the Perona-Malik equation \fer{eq:pm}. Then, for the cross-diffusion problem and for the gradient based Perona-Malik equation, that is, when $g\equiv g(\abs{\grad u})$, we use the finite element method (FEM) to discretize in space, with linear-wise basis. For the Laplacian based Perona-Malik equation, the use of FEM  would require to consider a basis of, at least,  piece-wise quadratic polynomials to approximate the second order derivatives. Instead of incrementing the order of the FEM basis, we preferred to use a simpler scheme for this filter, based on the finite differences method. 

We give the implementation details for the cross-diffusion system, being that of the Perona-Malik equation similar. Following \cite{Gilboa2001}, we fix the initial data as $\bu_0 = (u_{\rm noisy}, 0)$ and consider the edge detector function as given by 
\begin{equation*}
 g(\bu) = g(u_2),\qtext{with } g(s) = \exp\big(-s^2/\lambda_{cd}^2\big),
\end{equation*}
with $\lambda_{cd}>0$ optimized according to each method and image. The right hand side of \fer{eq:pde} is taken as $f_i\equiv 0$ (pure filtering case). 
Finally, like in \cite{Gilboa2001}, the diffusion coefficient matrix is taken as, for fixed $\theta \in [0,\pi/2)$,
\[
 A = \begin{pmatrix} \cos \theta & -\sin \theta  \\ \sin \theta  & \cos \theta  \end{pmatrix},
\]
with $\theta = \pi/30$ in the experiments.

We used the open source software \texttt{deal.II} \cite{dealII83} to implement a time semi-implicit 
scheme with spatial linear-wise finite element discretization. For the time discretization, we take 
in the experiments a uniform partition of $[0,T]$ of time step $\tau$. For the spatial  discretization, we take the natural uniform partition of the rectangle $[0,W]\times[0,H]$, where $W$ and $H$ stand for the width and height (in pixels) of the given image,  respectively.

Let, initially, $t=t_0=0$ and  set $u_{i}^0=u_{i0}$. For $n\geq 1$ the problem is: Find $u_{ i}^{n}\in S^h$ such that for $i=1,2$,  
\begin{equation}\label{eq:pde_discr.s4}
\begin{array}{l}
\frac{1}{\tau}\big( u^n_{ i}-u^{n-1}_{ i} , \chi )^h
+ \big( g(u_2^n)(a_{i1} \grad u_1^n + a_{i2}\grad u_2^n  ) ,\grad\chi \big)^h =0,
\end{array}
\end{equation}
for every $ \chi\in S^h $, the finite element space of piecewise $\mathbb{Q}_1$-elements. 
Here, $(\cdot,\cdot)^h$ 
stands for a discrete semi-inner product on $\mathcal{C}(\overline{\Omega} )$.

Since \fer{eq:pde_discr.s4} is a nonlinear algebraic problem, we use a fixed point argument to 
approximate its solution,  $(u_{1}^n,u_{2}^n)$, at each time slice $t=t_n$, from the previous
approximation $u_{ i}^{n-1}$.  Let $u_{ i}^{n,0}=u_{i}^{n-1}$. 
Then, for $k\geq 1$ the linear problem to solve is: Find $u_{ i}^{n,k}$ such that for 
$i=1,2$, and for all $\chi \in S^h$ 
\begin{equation*}
\begin{array}{l}
 \frac{1}{\tau}\big( u^{n,k}_{ i}-u^{n-1}_{  i} , \chi )^h
+ \big( g(u_2^{n,k-1})(a_{i1} \grad u_1^{n,k} + a_{i2}\grad u_2^{n,k}  ) ,\grad\chi \big)^h =0.
\end{array}
\end{equation*}
We use the stopping criteria 
 \begin{equation}
 \label{tol}
 \max _{i=1,2} \nor{u_{ i}^{n,k}-u_{ i}^{n,k-1}}_2 <\text{tol},
 \end{equation}
for values of $\text{tol}$ chosen empirically, and set $u_i^n=u_i^{n,k}$. See Algorithm~\ref{al.cd} for implementation details.

Turning to the neighborhood one-step filters, we used the simplest version of the Bilateral filter, also known as Yaroslavsky filter \cite{Yaroslavsky1985}, given by 
\begin{align*}
 \Ya u(\bx) = \frac{1}{C(\bx)} \int_{B_\rho(\bx)} \exp\Big(- \frac{\abs{u(\bx)-u(\by)}^2}{h^2} \Big) u(\by) d \by,   
\end{align*}
where $B_\rho(\bx)$ is a box of diameter $4\rho$ centered at $\bx$, and $\rho,~h$ are positive parameters. The term $C(\bx)$ is the normalizing factor
\begin{align*}
 C(\bx) = \int_{B_\rho(\bx)} \exp\Big(- \frac{\abs{u(\bx)-u(\by)}^2}{h^2} \Big) d \by.   
\end{align*}
For its implementation, we used the fast algorithm introduced in \cite{Porikli2008, Galiano2015, Galiano2015b}.

The Nonlocal Means filter is defined as 
\begin{align*}
 \NLM u(\bx) = \frac{1}{C(\bx)} 
 \int_\O \exp\Big(-\frac{M_\sigma(\bx,\by)}{h^2} \Big) u(\by) d \by,   
\end{align*}
with $h>0$, and the normalizing factor
\begin{align*}
C(\bx) =  \int_\O \exp\Big(-\frac{M_\sigma(\bx,\by)}{h^2} \Big)d \by.  
\end{align*}
The nonlocal term, $M_\sigma$, is given in convolution form by
\begin{align*}
M_\sigma(\bx,\by) = \int_{\R^2} G_\sigma(\bz) \abs{u(\bx+\bz)-u(\by+\bz)}^2 d\bz,
\end{align*}
where $G_\sigma$ is a Gaussian of standard deviation $\sigma$. In practice, the parameter $h$ is fixed in terms of $\sigma$, so this method has a single parameter to be optimized.  For the discrete implementation, we used the patch-wise approach introduced by the authors \cite{Buades2011}.

\begin{algorithm}
\caption{FEM-Fixed point method for cross-diffusion problem \fer{eq:pde_discr.s4}  }\label{al.cd}
\begin{algorithmic}[1]
\Require{Noisy image: $u_{noisy}$; Final time: $T$; Time step: $\tau$}
\Require{Diffusion matrix: $a$; Edge detector function: $g$}
\Require{Tolerance $tol$ for the fixed point method; Max. number of fixed point iterations: maxFpIter}
\Require{Spatial domain}
\State{Compute triangulation, function basis $\vfi$, and massMatrix}
\State{$\bu = (u_{noisy}, 0)$, $t=0$}
\While {$t<T$}
  \State{$t ~ +\!=  ~\tau$}
  \State{$\bu_{old} = \bu$, $\bu_{fp} = \bu$, relDifference = $\infty$, fpIter=0}
  \While{relDifference $>$ tol \textbf{and}  fpIter $<$ maxFpIter}
    \State{fpIter $+\!=$ 1}
    \State{$\bu_{fp,old} = \bu_{fp}$}
    \For{All cells of the triangulation}
      \State{cellMatrix = 0}
      \For{All points in each cell} 
        \State{$edge = g(u_{fp,old, 2})$} 
        \For{$\vfi_i$'s of cell basis}
          \For{$\vfi_j$'s of cell basis}
            \State{$cellMatrix(i,j) +\!= \tau*edge*a(i,j)*\grad \vfi_i *\grad \vfi_j * area$}
          \EndFor
        \EndFor
      \EndFor
      \State{Add cellMatrix contribution to systemMatrix}
    \EndFor
    \State{systemMatrix $+\!=$ massMatrix  }  
    \If{fpIter = 1}
      \State{RHS = systemMatrix*$\bu$}
    \EndIf
    \State{Solve by LU method:  systemMatrix*$\bu_{fp}$ = RHS}
    \State{Compute relative difference between $\bu_{fp}$ and $\bu_{fp,old}$}
  \EndWhile  
  \State{ $\bu = \bu_{fp}$}
\EndWhile
\Ensure{ $\bu$}
\end{algorithmic}
\end{algorithm}

\subsection{Experiment data and results}

\begin{table}[h!]
{\footnotesize
\centering
\begin{tabular}{|c||c|c|c|c|c|c|}
\hline
       & Initial & CD & PM-L & PM-G & BF& NLM \\
\hline 
\hline 
Parameters  & & $T$, $\lambda_{cd}$ &  $T$, $\lambda_{pm}$ & $T$, $\lambda_{pm}$ & $h$, $\rho$ & $\sigma$ \\
 \hline 
\multicolumn{2}{|c}{} & \multicolumn{3}{c}{\emph{added value} } & \multicolumn{2}{c|}{} \\
\hline 
Opt. Par.&  & (0.2, 0.15) & (0.8, 10) & (0.3, 20) &  (64, 4) & 8 \tabularnewline \hline 
PSNR & 31.8618 & 35.5620 & 34.0086 & 36.7013 & 36.2672 & 38.4068 \tabularnewline \hline 
NCC & 0.9994 & 0.99975 & 0.99963 & 0.9998 & 0.99978 & 0.99987 \tabularnewline \hline 
SSIM & 0.99503 & 0.99788 & 0.99695 & 0.99838 & 0.99819 & 0.99891 \tabularnewline \hline 
\multicolumn{2}{|c}{} & \multicolumn{3}{c}{\emph{boat} } & \multicolumn{2}{c|}{} \\
\hline 
Opt. Par.&  & (0.1, 0.1) & (0.3, 10) & (0.1, 20) &  (12, 3) & 5 \tabularnewline \hline 
PSNR & 34.7000 & 36.2916 & 35.4629 & 36.4356 & 36.0391 & 36.9337 \tabularnewline \hline 
NCC & 0.99942 & 0.9996 & 0.99952 & 0.99962 & 0.99957 & 0.99966 \tabularnewline \hline 
SSIM & 0.99502 & 0.99653 & 0.99579 & 0.99666 & 0.9963 & 0.99702 \tabularnewline \hline 
\multicolumn{2}{|c}{} & \multicolumn{3}{c}{\emph{clock} } & \multicolumn{2}{c|}{} \\
\hline 
Opt. Par.&  & (0.14, 0.15) & (0.5, 10) & (0.3, 20) &  (33, 4) & 7 \tabularnewline \hline 
PSNR & 32.9488 & 35.5378 & 34.187 & 37.0526 & 36.8077 & 38.4792 \tabularnewline \hline 
NCC & 0.99957 & 0.99976 & 0.99968 & 0.99983 & 0.99982 & 0.99988 \tabularnewline \hline 
SSIM & 0.99505 & 0.99726 & 0.99625 & 0.99809 & 0.99796 & 0.99863 \tabularnewline \hline 
\multicolumn{2}{|c}{} & \multicolumn{3}{c}{\emph{house} } & \multicolumn{2}{c|}{} \\
\hline 
Opt. Par.&  & (0.13 , 0.15) & (0.5, 10) & (0.2, 20) &  (20, 3) & 5 \tabularnewline \hline 
PSNR & 34.4251 & 37.0206 & 35.7168 & 37.8745 & 37.7222 & 39.2158 \tabularnewline \hline 
NCC & 0.99859 & 0.99924 & 0.99891 & 0.99937 & 0.99934 & 0.99954 \tabularnewline \hline 
SSIM & 0.99504 & 0.99726 & 0.99608 & 0.99777 & 0.99769 & 0.99838 \tabularnewline \hline 
\multicolumn{2}{|c}{} & \multicolumn{3}{c}{\emph{test} } & \multicolumn{2}{c|}{} \\
\hline 
Opt. Par.&  & (0.12, 0.1) & (0.3, 10) & (0.3, 25) &  (17, 64) & 14 \tabularnewline \hline 
PSNR & 28.8457 & 30.1614 & 29.1625 & 30.8796 & 32.5611 & 33.0273 \tabularnewline \hline 
NCC & 0.99874 & 0.99916 & 0.9973 & 0.99928 & 0.99947 & 0.99964 \tabularnewline \hline 
SSIM & 0.99711 & 0.99786 & 0.99408 & 0.9982 & 0.99881 & 0.9989 \tabularnewline \hline 
\multicolumn{2}{|c}{} & \multicolumn{3}{c}{\emph{tree} } & \multicolumn{2}{c|}{} \\
\hline 
Opt. Par.&  & (0.12, 0.15) & (0.4, 10) & (0.25, 20) &  (35, 4) & 8 \tabularnewline \hline 
PSNR & 31.7344 & 33.7960 & 32.5588 & 34.7020 & 34.1934 & 35.2368 \tabularnewline \hline 
NCC & 0.99795 & 0.99874 & 0.99617 & 0.99897 & 0.99884 & 0.9991 \tabularnewline \hline 
SSIM & 0.99517 & 0.99697 & 0.99092 & 0.99756 & 0.99727 & 0.99785 \tabularnewline \hline 
\end{tabular}
\caption{Optimal parameters and quality results for the set of natural images.} 
\label{table1}
}
\end{table}

\begin{table}[h!]
{\footnotesize
\centering
\begin{tabular}{|c|c||c|c||c|c||c|}
\hline
Measure  & Initial & CD & PM-L & PM-G & BF& NLM \\
\hline 
\hline 
Parameters  & & $T$, $\lambda_{cd}$ &  $T$, $\lambda_{pm}$ & $T$, $\lambda_{pm}$ & $h$, $\rho$ & $\sigma$ \\
 \hline 
\multicolumn{2}{|c}{} & \multicolumn{3}{c}{\emph{bark} } & \multicolumn{2}{c|}{} \\
\hline 
Opt. Par.&  & (0.02, 0.3) & (0.05, 50) & (0.03, 70) &  (12, 6) & 8 \tabularnewline \hline 
PSNR & 30.8550 & 31.3764 & 30.9596 & 31.3815 & 30.9778 & 30.8820 \tabularnewline \hline 
NCC & 0.99877 & 0.99892 & 0.99873 & 0.99892 & 0.99881 & 0.99878 \tabularnewline \hline 
SSIM & 0.99517 & 0.99566 & 0.99492 & 0.99568 & 0.99528 & 0.9952 \tabularnewline \hline 
\multicolumn{2}{|c}{} & \multicolumn{3}{c}{\emph{bricks} } & \multicolumn{2}{c|}{} \\
\hline 
Opt. Par.&  & (0.01, 0.25) & (0.05, 50) & (0.04, 50) &  (13, 5) & 4 \tabularnewline \hline 
PSNR & 30.8471 & 31.1076 & 30.975 & 31.1590 & 31.1448 & 30.9597 \tabularnewline \hline 
NCC & 0.99878 & 0.99886 & 0.99881 & 0.99887 & 0.99886 & 0.99881 \tabularnewline \hline 
SSIM & 0.99517 & 0.9954 & 0.99522 & 0.99547 & 0.99547 & 0.9953 \tabularnewline \hline 
\multicolumn{2}{|c}{} & \multicolumn{3}{c}{\emph{bubbles} } & \multicolumn{2}{c|}{} \\
\hline 
Opt. Par.&  & (0.03, 0.1) & (0.05, 30) & (0.04, 40) &  (7, 3) & 2 \tabularnewline \hline 
PSNR & 36.1726 & 36.8475 & 36.3222 & 37.0446 & 36.3846 & 36.2160 \tabularnewline \hline 
NCC & 0.99944 & 0.99953 & 0.99946 & 0.99955 & 0.99946 & 0.99945 \tabularnewline \hline 
SSIM & 0.99502 & 0.99569 & 0.99515 & 0.9959 & 0.99518 & 0.99507 \tabularnewline \hline 
\multicolumn{2}{|c}{} & \multicolumn{3}{c}{\emph{holes} } & \multicolumn{2}{c|}{} \\
\hline 
Opt. Par.&  & (0.03, 0.3) & (0.05, 20) & (0.15, 30) &  (14, 7) & 14 \tabularnewline \hline 
PSNR & 29.2774 & 29.9741 & 29.2695 & 29.9926 & 30.1212 & 30.1021 \tabularnewline \hline 
NCC & 0.999 & 0.9992 & 0.99892 & 0.99916 & 0.99914 & 0.9992 \tabularnewline \hline 
SSIM & 0.99617 & 0.99672 & 0.99584 & 0.99677 & 0.99681 & 0.99686 \tabularnewline \hline 
\multicolumn{2}{|c}{} & \multicolumn{3}{c}{\emph{straw} } & \multicolumn{2}{c|}{} \\
\hline 
Opt. Par.&  & (0.01, 0.3) & (0.01, 50) & (0.02, 70) &  (12, 5) & 14 \tabularnewline \hline 
PSNR & 30.8735 & 31.1492 & 30.8961 & 31.1175 & 30.9911 & 30.8918 \tabularnewline \hline 
NCC & 0.99878 & 0.99886 & 0.99875 & 0.99885 & 0.99881 & 0.99878 \tabularnewline \hline 
SSIM & 0.99518 & 0.99542 & 0.99506 & 0.99541 & 0.99528 & 0.9952 \tabularnewline \hline 
\end{tabular}
\caption{Optimal parameters and quality results for the set of texture images.} 
\label{table2}
}
\end{table}

For the cross-diffusion and the Perona-Malik equations, the time step is fixed as $\tau = 0.01$, and the tolerance for the fixed point iteration inside each time loop, see \fer{tol}, is taken  as ${\rm tol} = 10^{-3}$.

In Table~\ref{table1} we show the optimal parameters found for each method, in the sense of PSNR maximization, and the resulting quality measures for each image in the set of natural images. For most of these images, the Laplacian based Perona-Malik equation gives the poorest results. The cross-diffusion, the gradient based Perona-Malik equation and the Bilateral (Yaroslavsky) filters give, in general, similar results. The Nonlocal Means filter outperforms the other methods for all the images. 

In Table~\ref{table2} we present the same information for the set of texture images. 
In this case the results are not conclusive. The relative performance among methods is closer than in the case of natural images and, in fact, the cross-diffusion and the gradient based Perona-Malik equations give the best results for most of the images. 
However, the denoising capacity of all of them is very limited and the denoising effect is hard to visualize. 

\begin{figure*}[ht]
\centering
 {\includegraphics[width=0.30\textwidth,height=0.25\textheight]{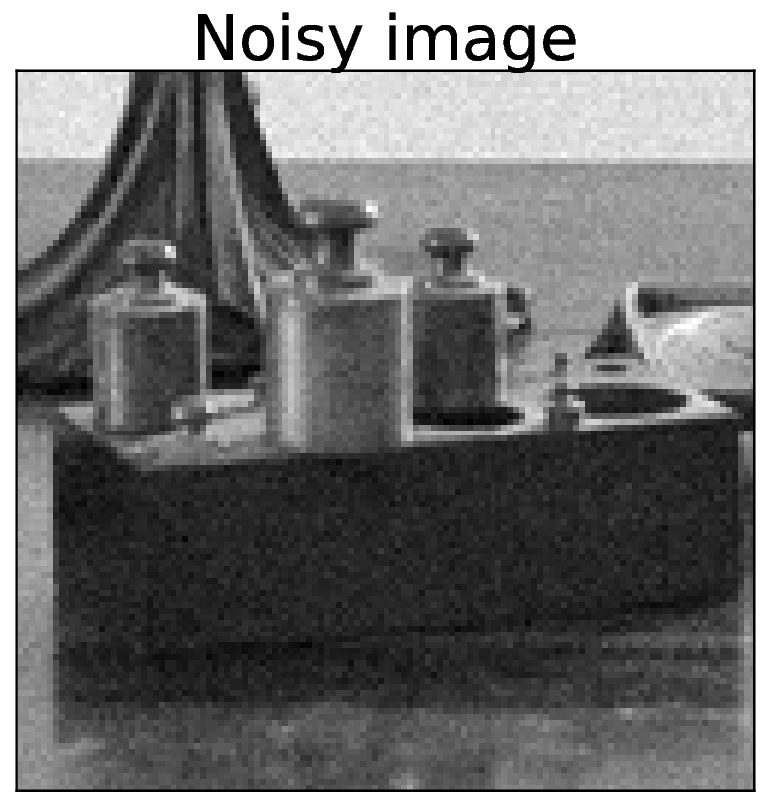}}
 {\includegraphics[width=0.30\textwidth,height=0.25\textheight]{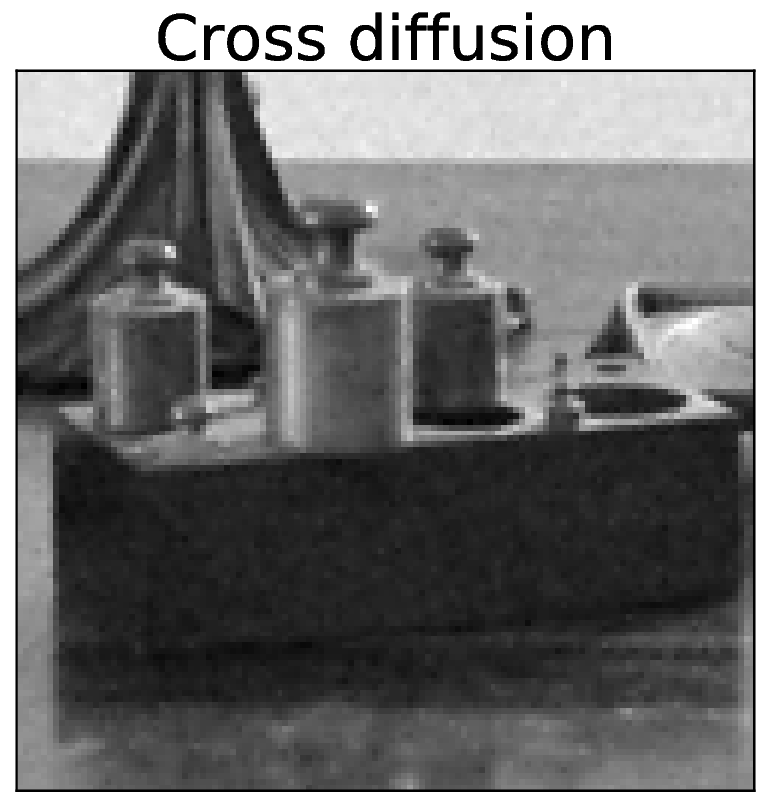}}
 {\includegraphics[width=0.30\textwidth,height=0.25\textheight]{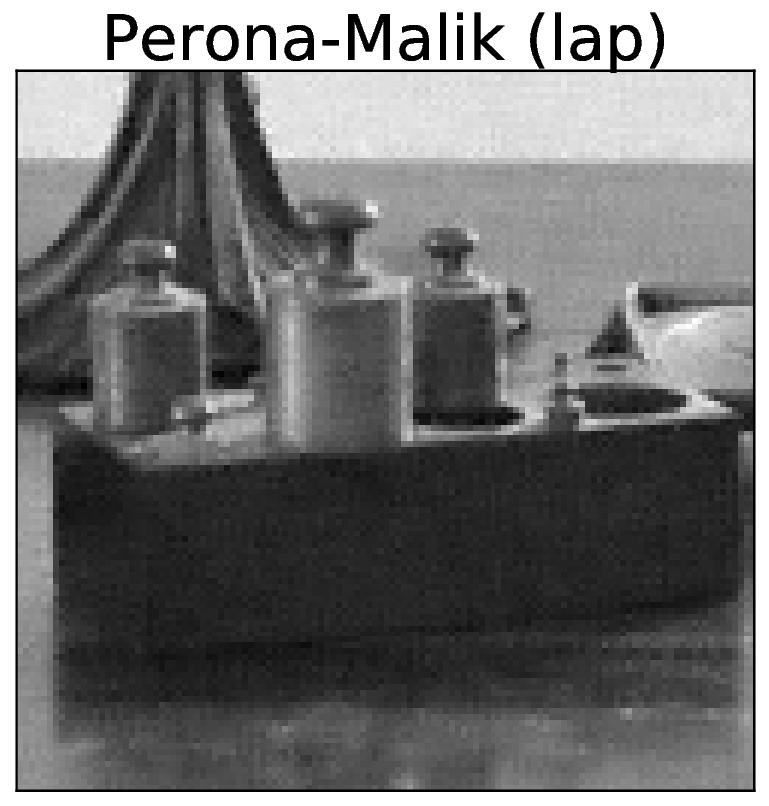}}\\
 {\includegraphics[width=0.30\textwidth,height=0.25\textheight]{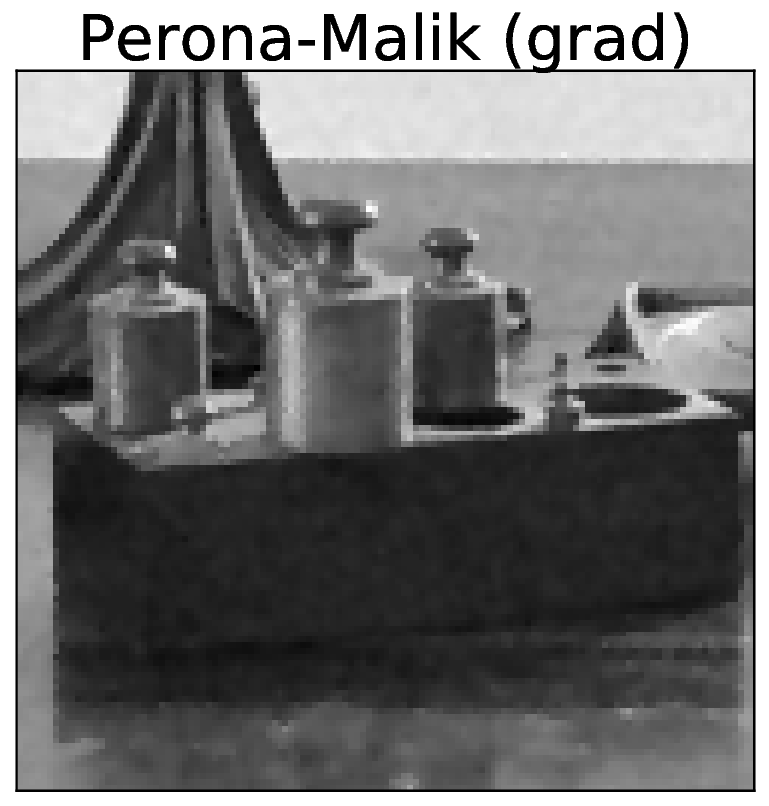}}
 {\includegraphics[width=0.30\textwidth,height=0.25\textheight]{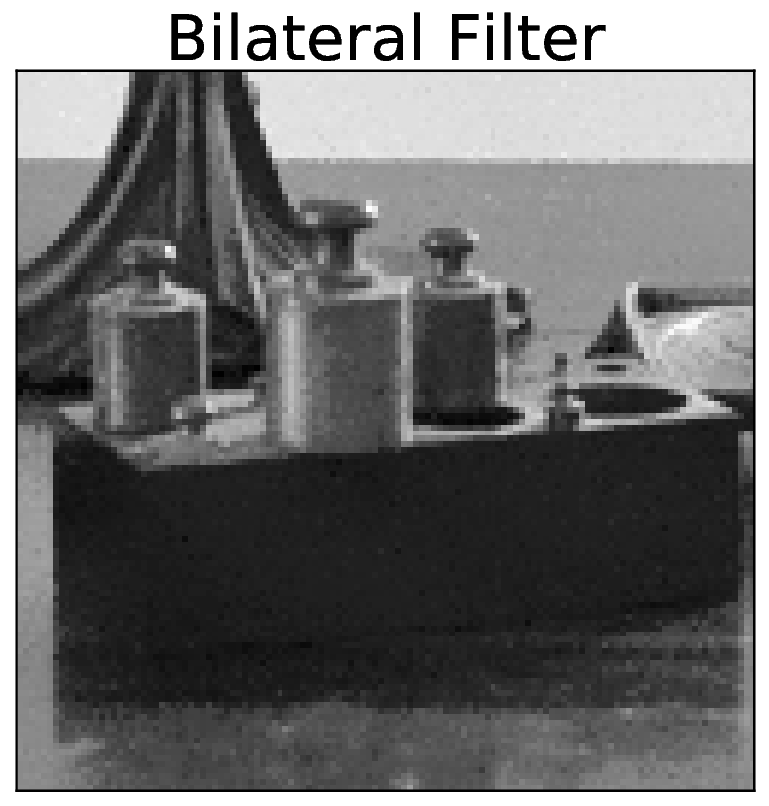}}
 {\includegraphics[width=0.30\textwidth,height=0.25\textheight]{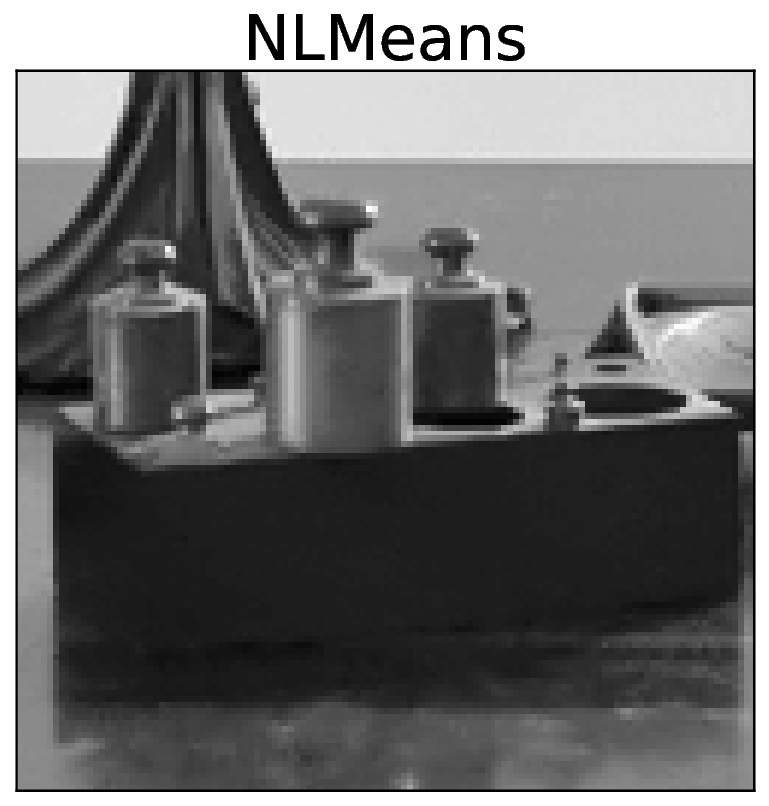}}
   \caption{{Detail of the natural image \emph{added value} showing the performance of each  method. }} 
\label{fig2}
\end{figure*}

\begin{figure*}[ht]
\centering
 {\includegraphics[width=0.30\textwidth,height=0.25\textheight]{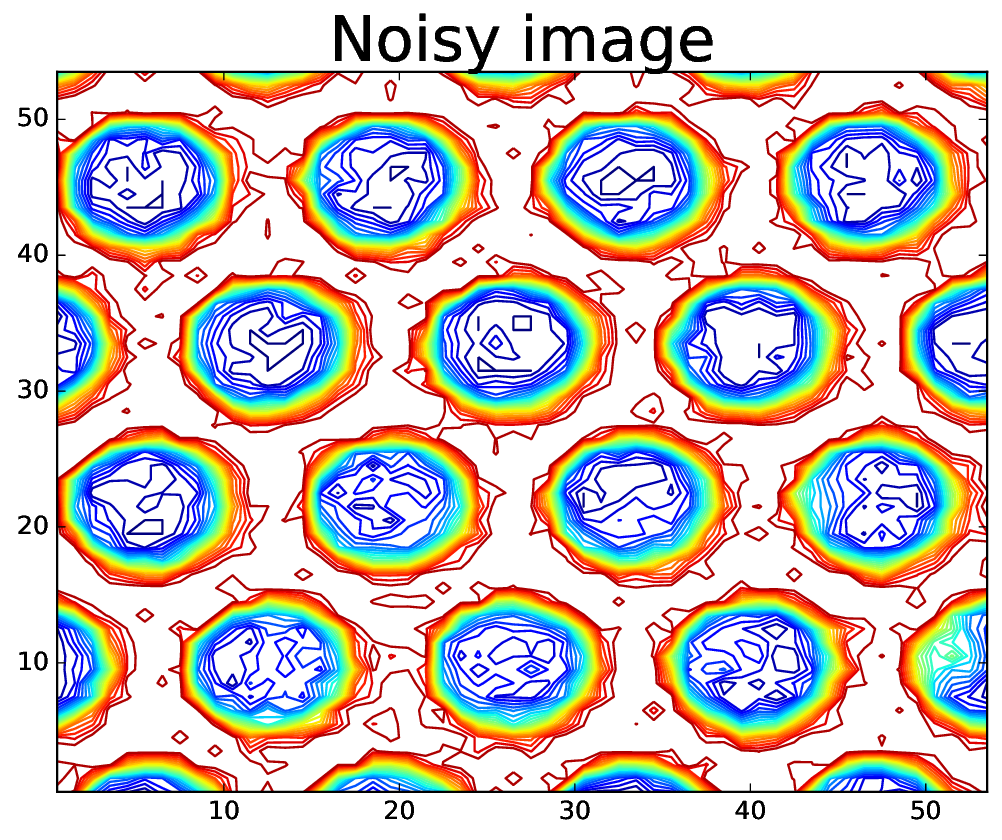}}
 {\includegraphics[width=0.30\textwidth,height=0.25\textheight]{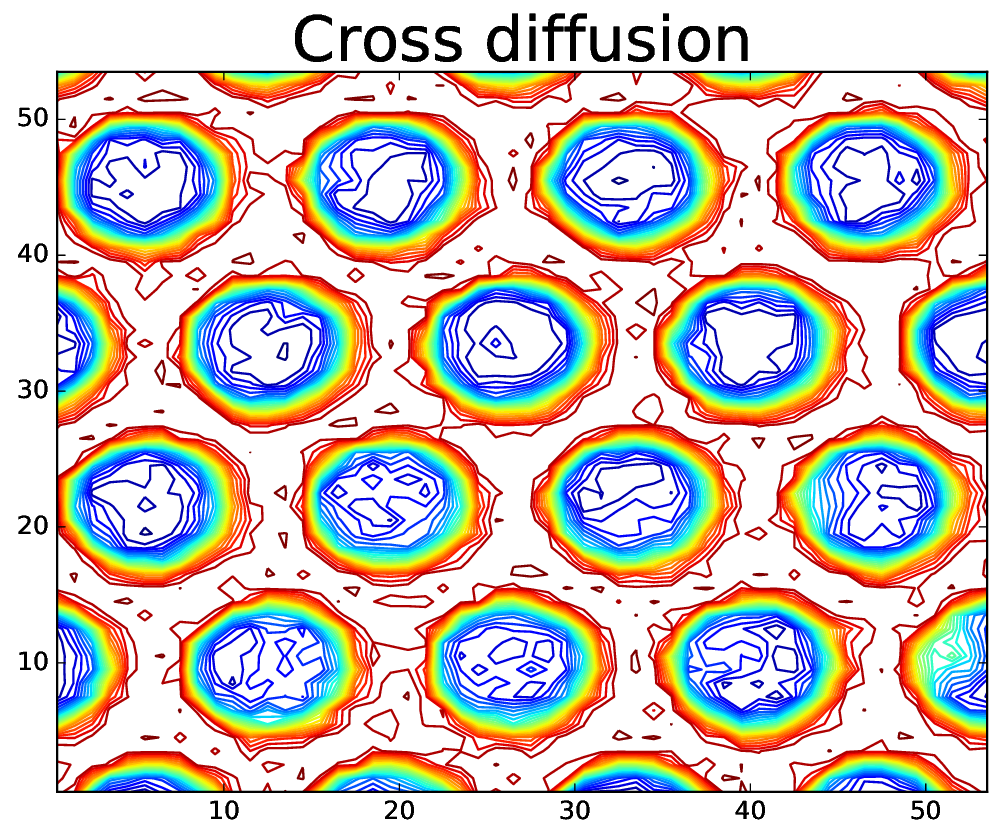}}
 {\includegraphics[width=0.30\textwidth,height=0.25\textheight]{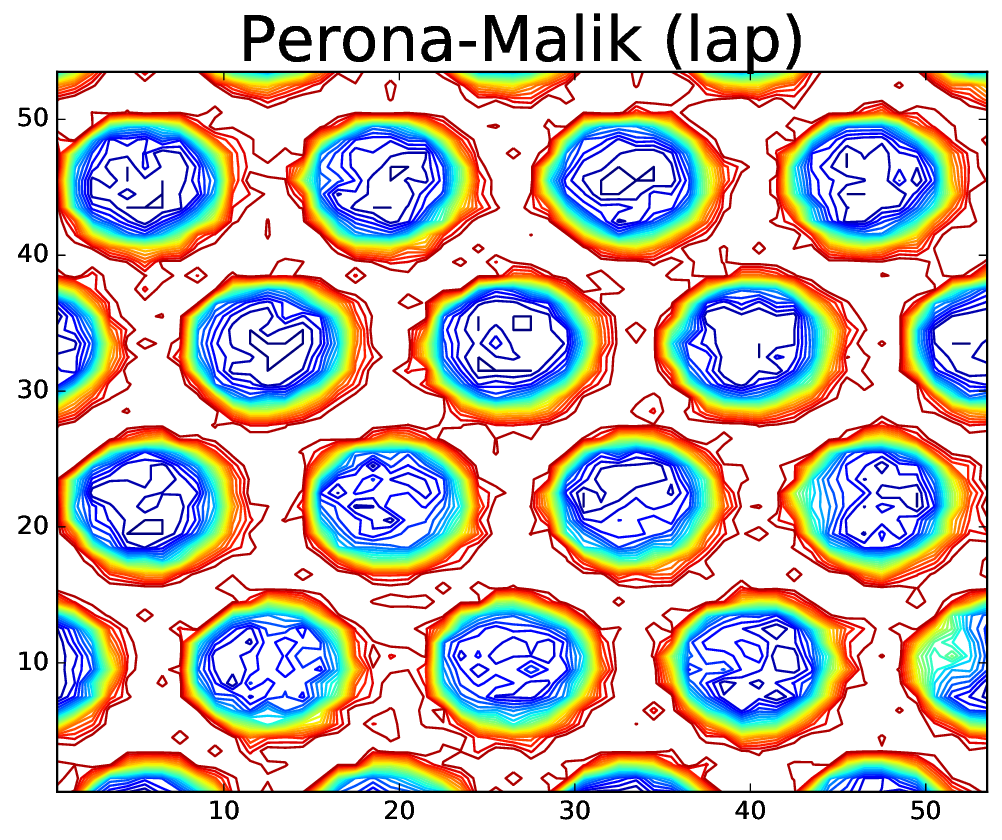}}\\
 {\includegraphics[width=0.30\textwidth,height=0.25\textheight]{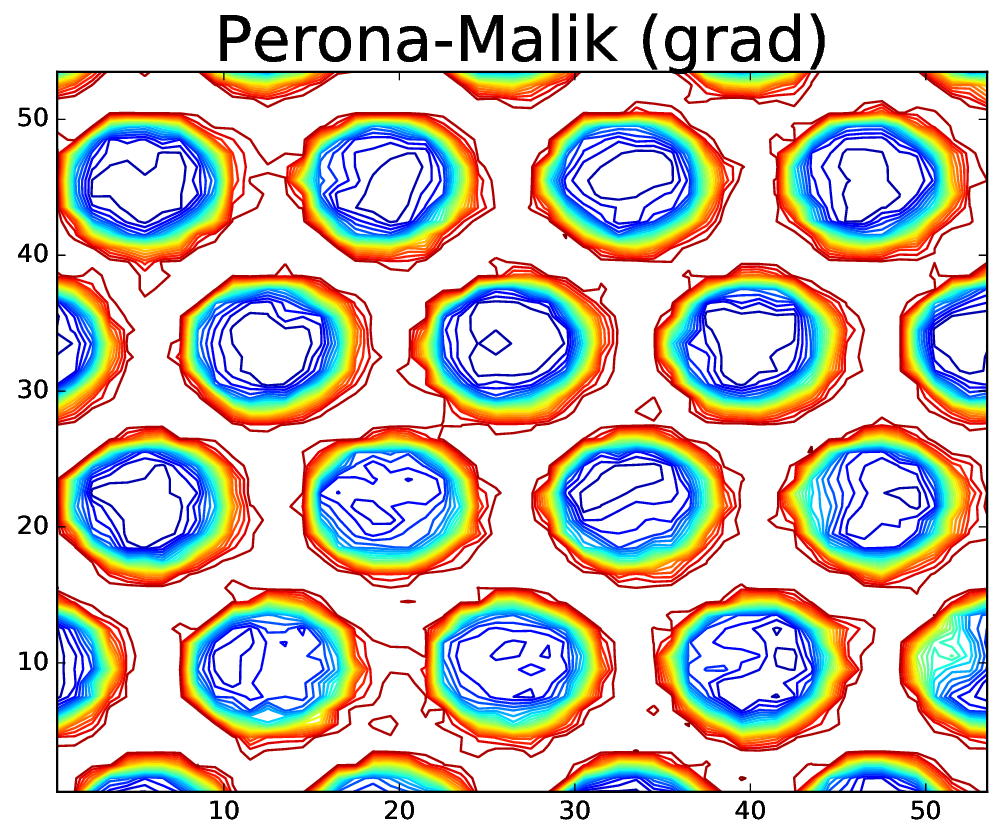}}
 {\includegraphics[width=0.30\textwidth,height=0.25\textheight]{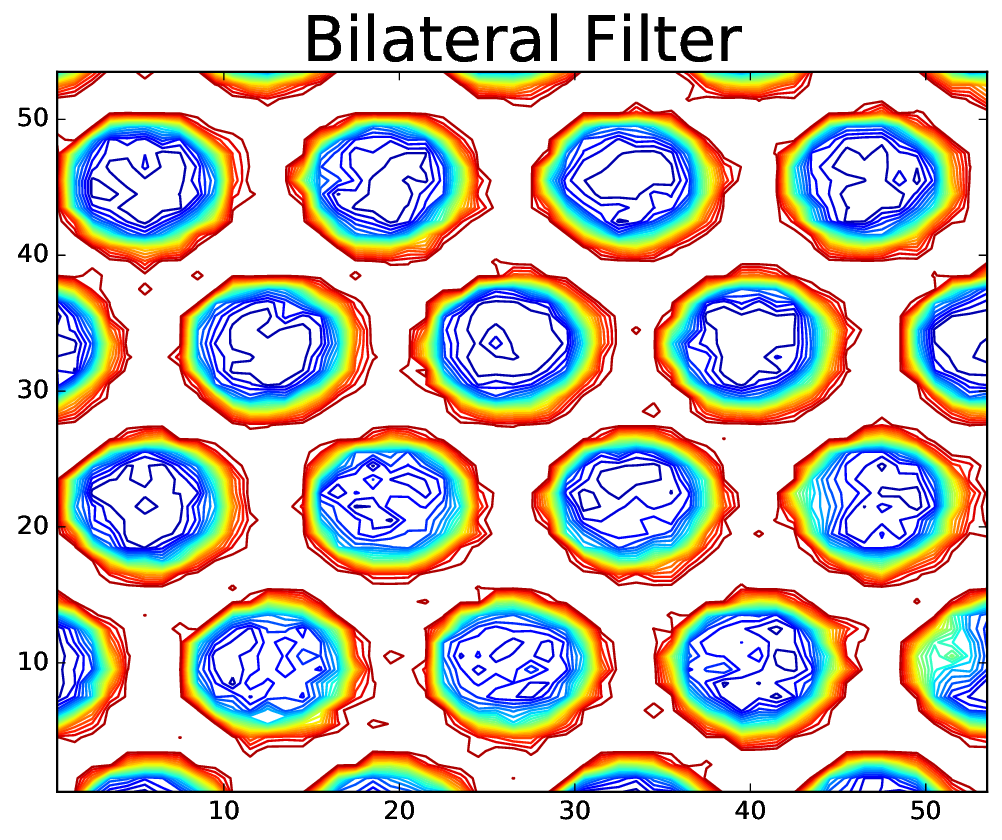}}
 {\includegraphics[width=0.30\textwidth,height=0.25\textheight]{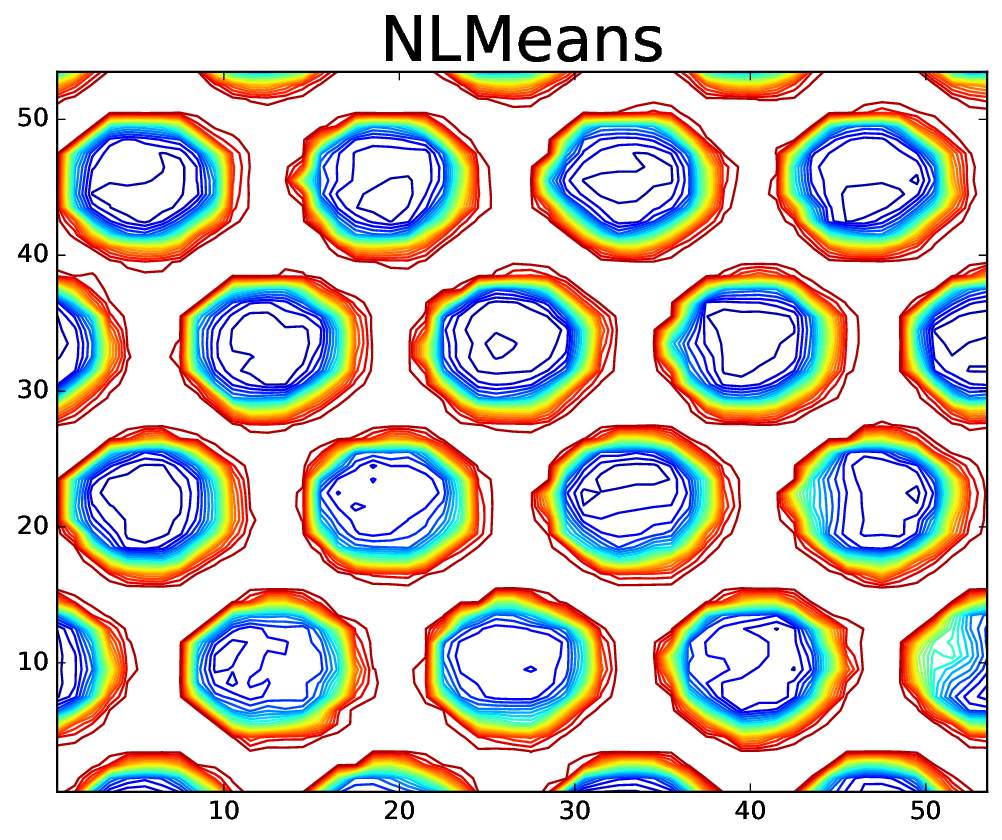}}
   \caption{Contours plots for a detail of the texture image \emph{holes}.} 
\label{fig3}
\end{figure*}

In Figure~\ref{fig2} we show a detail of the image \emph{added value} showing the performance of all the methods for a natural image.  In Figure~\ref{fig3}
we show the contours plots for a detail of the texture image \emph{holes}.

\section{Proofs}

\no\textbf{Proof of Theorem~\ref{th:steady}. }   
We start considering the following problem: Find $\bu:\O\to\R^2$ such that, for $i=1,2$, and $i\neq j$,   
\begin{align}
 &\gamma_{i} u_i  -  \frac{\det(A)}{a_{jj}}\Div \big( g( \bu) \grad u_i\big) 
 = G_i-\alpha_{ij}G_j +\alpha_{ij} \gamma_{j} u_j && \qtext{in }\O, 
  \label{eq:pde:aux:1}\\
& g( \bu) \grad u_i \cdot \nu =0 && \qtext{on }\partial\O.
\label{eq:bc.aux:1}
\end{align}
with $\alpha_{ij} = a_{ij}/a_{jj}$. Notice that, due to assumption (H)$_2$,
\begin{align*}
\det (A) >0 ,\quad a_{ii}>0,\qtext{and}\quad \abs{\alpha_{ij}} < 1.
\end{align*}

\smallskip

\no\emph{Step 1: Linearizing.}
Let $\tilde \bu \in L^2(\O)^2$ be given and consider the function $g_k(s_1,s_2) = g(T_k(s_1),T_k(s_2))$, with 
the truncation function $T_k:\R\to\R$ given by 
\begin{align*}
 T_{k}(s) = \left\{\begin{array}{ll}
                      -k & \text{if } s\leq -k, \\
                      s & \text{if } -k\leq s\leq k, \\
                      k & \text{if } s\geq k.
                     \end{array}
\right.
\end{align*}
Let $\sigma \in [0,1]$. We set the  following linear problem: Find $\bu \in H^1(\O)^2$ such that
for $i=1,2$, and $i\neq j$,   
\begin{align}
 &\gamma_{i} u_i  -  \frac{\det(A)}{a_{jj}}\Div \big( g_k( \tilde \bu) \grad u_i\big) 
 = \sigma(G_i-\alpha_{ij}G_j) +\alpha_{ij} \gamma_{j} u_j && \qtext{in }\O, 
  \label{eq:pde:aux:2}\\
& g_k( \tilde \bu) \grad u_i \cdot \nu =0 && \qtext{on }\partial\O.
\label{eq:bc:aux:2}
\end{align}
Since, $g_k(\R^2)\subset g(K)$, with $K= [-k,k]\times [-k,k]$, and by assumption  (H)$_3$, $g(K)\subset (0,\infty)$, that is, the problem is uniformly elliptic, we may apply \cite{LU}[Chapter 6, Theorem 2.1], to deduce the existence of a unique weak solution $\bu \in H^1(\O)^2 \cap L^\infty(\O)^2$ of problem \fer{eq:pde:aux:2}-\fer{eq:bc:aux:2}.

We now show that the $L^\infty(\O)$ bound does not depend on $\tilde \bu$ or $k$.
Let $[s]_+$ and $[s]_-$ denote the positive and negative parts of $s$, so that 
$s=[s]_+ - [s]_-$. 
Taking, for $p\geq 1$,  $[u_i]_+ ^p \in H^1(\O)$ as test function in the weak formulation of \fer{eq:pde:aux:2}-\fer{eq:bc:aux:2}, and using Young's inequlity, we deduce 
\begin{align*}
 \gamma_{i} \nor{[u_i]_+}_{p+1} \leq   \sigma \big(\nor{G_i}_{p+1} + \abs{\alpha_{ij}} \nor{ G_j}_{p+1} \big)+  \gamma_j \abs{\alpha_{ij}} \nor{[u_j]_+}_{p+1} .
\end{align*}
Summing for $i,j=1,2$, with $i\neq j$, we obtain 
\begin{align*}
 \gamma_{1}(1-\abs{\alpha_{21}}) & \nor{[u_1]_+}_{p+1}+
  \gamma_{2}(1-\abs{\alpha_{12}})\nor{[u_2]_+}_{p+1} \nonumber \\
  &\leq   \sigma \big(1+\abs{\alpha_{21}}\big) \nor{G_1}_{p+1} + 
  \sigma\big(1+\abs{\alpha_{12}}\big) \nor{ G_2}_{p+1} 
\end{align*}
A similar estimate holds for $[u_i]_- ^{p+1}$, so we obtain, taking the limit $p\to\infty$,
\begin{align}
 \gamma_{1}(1-\abs{\alpha_{21}})\nor{u_1}_{\infty}+ &
  \gamma_{2}(1-\abs{\alpha_{12}})\nor{u_2}_{\infty} \nonumber \\
&   \leq   \sigma\big(1+\abs{\alpha_{21}}\big) \nor{G_1}_{\infty} + 
  \sigma\big(1+\abs{\alpha_{12}}\big) \nor{ G_2}_{\infty} .
 \label{ineq:linf:init}
\end{align}

 \smallskip

\no\emph{Step 2: Fixed point.} Consider the operator $P:L^2(\O)^2 \times[0,1]\to L^2(\O)^2$ given by $P(\tilde \bu, \sigma ) = \bu$, being $\bu$ the solution of the linear problem \fer{eq:pde:aux:2}-\fer{eq:bc:aux:2} corresponding to $(\tilde \bu,\sigma)$. We check that $P$ satisfies the conditions of the Leray-Schauder's Fixed Point theorem, that is 
\begin{enumerate}
 \item $P$ is continuous and compact.
 
 \item $P(\bv,0)=0$ for all $\bv\in L^2(\O)^2$.
 
 \item The fixed points of $P$ are uniformly bounded in $ L^2(\O)^2$.
\end{enumerate}
1. Let $\bv_n \in L^2(\O)$ be a sequence strongly convergent to $\bv\in L^2(\O)$, and 
$\sigma_n \to \sigma$. Define $\bu_n = P(\bv_n,\sigma_n)$. Using $u_{i,n}\in H^1(\O)$ as a test function in the weak formulation of \fer{eq:pde:aux:2}-\fer{eq:bc:aux:2} leads to 
\begin{align}
 \gamma_i \int_\O \abs{u_{i,n}}^2 + \frac{\det(A)}{a_{jj}} \min(g(K)) \int_\O \abs{\grad u_{i,n}}^2  \leq & \sigma_n \int_\O (\abs{G_i} + \abs{\alpha_{ij}} \abs{G_j} )\abs{u_{i,n}} \nonumber \\ 
 & + \abs{\alpha_{ij}} \int_\O \abs{u_{i,n}}\abs{u_{j,n}}. \label{aqui}
\end{align}
We first use H\"older's inequality to deduce, after summing for $i,j=1,2$, with $i\neq j$
\begin{align*}
 \gamma_{1}(1-\abs{\alpha_{21}}) & \nor{u_1}_{2}^2+
  \gamma_{2}(1-\abs{\alpha_{12}})\nor{u_2}_{2}^2 \nonumber \\
  &\leq   \sigma_n \big(1+\abs{\alpha_{21}}\big) \nor{G_1}_{2}^2 + 
  \sigma_n \big(1+\abs{\alpha_{12}}\big) \nor{ G_2}_{2}^2 ,
\end{align*}
and then, from \fer{aqui}, we also deduce that $\nor{\grad u_{i,n}}_2$ is uniformly bounded. Therefore, there exists a subsequence of $\bu_n$ (not relabeled), and a function $\bu\in H^1(\O)^2$, such that 
\begin{align*}
& \grad u_{i,n} \wto \grad u_i &&\text{weakly in }L^2(\O),\\
& u_{i,n}\to u_i &&\text{strongly in }L^2(\O).
\end{align*}
 Since, by assumption, $v_{i,n}\to v_i$ in $L^2(\O)$ (and a.e. in $\O$, for a subsequence), and $\nor{g_k(v_{i,n})}_\infty$ is uniformly bounded with respect to $n$, 
 we deduce by the dominated convergence theorem that $g_k(v_{i,n})\to g_k(v_i)$ in $L^q(\O)$, for all $q<\infty$. Therefore,
\begin{align*}
 g_k(\bv_n) \grad u_{i,n} \wto g_k(\bv) \grad u_i \qtext{weakly in } L^\gamma(\O)
\end{align*}
for any $\gamma <2$. But since the sequence $g_k(\bv_n) \grad u_{i,n}$ is bounded in $L^2(\O)$, we deduce that, in fact, 
\begin{align*}
 g_k(\bv_n) \grad u_{i,n} \wto g_k(\bv) \grad u_i \qtext{weakly in } L^2(\O).
\end{align*}
The other terms in the weak formulation of \fer{eq:pde:aux:2}-\fer{eq:bc:aux:2} (with the replacements $\tilde \bu = \bv_n$, $\bu = \bu_n$, and $\sigma=\sigma_n$) pass clearly to their corresponding limits. 
Finally, the uniqueness  of solutions of problem \fer{eq:pde:aux:2}-\fer{eq:bc:aux:2} implies that the whole sequence $\bu_n$ converges to $\bu$. So the continuity follows.
The compactness is directly deduced from the compact embedding $H^1(\O)\subset L^2(\O)$.

\smallskip

\no 2. If $\bu = P(\bv,0)$ then, in particular, $\bu$ satisfies estimate \fer{ineq:linf:init} with $\sigma=0$. Thus, $\bu = 0$ a.e. in $\O$. 

\smallskip

\no 3. If $\bu = S(\bu,\sigma)$, we again use that $\bu$ satisfies \fer{ineq:linf:init} with $\sigma \in [0,1]$, to get the uniform bound of $\bu$ in $L^2(\O)$.

\smallskip

Thus, the Leray-Schauder's fixed point theorem ensures the existence of a fixed point, $\bu$, of $P(\bu,1)$, which is a solution of \fer{eq:pde:aux:1}-\fer{eq:bc.aux:1} with $g$ replaced by $g_k$. But since this fixed point satisfies the $L^\infty(\O)$ bound 
\fer{ineq:linf:init} (with $\sigma=1$), we deduce that for $k$ large enough $g_k(\bu) = g(\bu)$, and thus $\bu$ is a solution of \fer{eq:pde:aux:1}-\fer{eq:bc.aux:1}.

\no\emph{Step 3: Solution of the original problem \fer{eqs:pde}-\fer{eqs:bc}. }
The weak solution, $\bu$, of problem \fer{eq:pde:aux:1}-\fer{eq:bc.aux:1} satisfies, for all $\vfi\in H^1(\O)$, 
\begin{align*}
 \gamma_1 \int_\O u_1 \vfi + \frac{\det(A)}{a_{22}} \int_\O g(\bu) \grad u_1 \cdot \grad \vfi =  \int_\O G_1 \vfi  - \alpha_{12} \int_\O G_2 \vfi + 
  \alpha_{12} \gamma_2 \int_\O u_2 \vfi,\\
 \gamma_2 \int_\O u_2 \vfi + \frac{\det(A)}{a_{11}} \int_\O g(\bu) \grad u_2 \cdot \grad \vfi =  \int_\O G_2 \vfi  - \alpha_{21} \int_\O G_1 \vfi + 
  \alpha_{21} \gamma_1 \int_\O u_1 \vfi. 
\end{align*}
Multiplying the first equation by $a_{11}a_{22}/\det(A)$, the second by $a_{12}a_{11}/\det(A)$ and adding the results leads to 
\begin{align*}
 \gamma_1 \int_\O u_1 \vfi + \int_\O g(\bu) (a_{11}\grad u_1 + a_{12}\grad u_ 2)\cdot \grad \vfi =  \int_\O G_1 \vfi .
\end{align*}
A similar combination, also gives 
\begin{align*}
 \gamma_2 \int_\O u_2 \vfi + \int_\O g(\bu) (a_{21}\grad u_1 + a_{22}\grad u_ 2)\cdot \grad \vfi =  \int_\O G_2 \vfi ,
\end{align*}
so $\bu$ is a weak solution of \fer{eqs:pde}-\fer{eqs:bc}. $\Box$

\bigskip

 \no\textbf{Proof of Corollary~\ref{cor:qssa}.}  
 The proof follows easily from Theorem~\ref{th:steady} by taking 
 $\gamma_i = \beta_i + 1/\tau >0$ and, recursively, 
 $G_i = \beta_{i0}u_{i0}+ u_i^n /\tau \in L^\infty(\O)$. $\Box$

 \begin{remark}
The $L^\infty$ bound \fer{ineq:linf:init} of the steady state problem translates, in the case of the QSS approximation, to 
\begin{align*}
 \nor{u_1^{n+1}}_{\infty} + \nor{u_2^{n+1}}_{\infty}
  \leq C_1 \big(\nor{ u_1^{n}}_{\infty} + \nor{ u_2^{n}}_{\infty}\big)
      +C_2 \big(\nor{ u_{10}}_{\infty} + \nor{ u_{20}}_{\infty}\big),
\end{align*}
 with $C_1>1$ unless $\alpha_{12}=\alpha_{21}= 0$, that is, $a_{12}=a_{21}=0$. Thus, 
  this bound is not useful for passing to the limit $\tau\to0$, and thus to prove the existence of the evolution problem, except when $A$ is diagonal. 
 \end{remark}

\bigskip
 
\no\textbf{Proof of Corollary~\ref{cor:diagonal}. }

\no\emph{Case 1: $a_{12}=a_{21}=0$ (diagonal case).} Consider the time discretization 
\fer{eq:pde:ss}-\fer{eq:bc:ss}, for which the existence of a solution $\bu^{n+1}$ is guaranteed by 
Corollary~\ref{cor:qssa}. In this case, a finer $L^\infty$ estimate of each time slice may be 
obtained. Indeed, the 
$L^{p+1}$ estimate \fer{ineq:linf:init} of Step 1 of the proof of Corollary~\ref{cor:qssa} reduces to 
\begin{align*}
 (1+\tau\beta_i) \nor{u_i^{n+1}}_{p+1}  \leq   \nor{ u_i^n}_{p+1}
  + \tau\beta_i \nor{u_{10}}_{p+1} , 
 \end{align*}
 and thus we get
\begin{align*}
  \nor{u_i^{n+1}}_{\infty}
  \leq   \frac{1}{1+\tau\beta_i} \nor{ u_i^n}_{\infty}+ \frac{\tau\beta_i}{1+\tau\beta_i} \nor{u_{i0}}_{\infty}.
 \end{align*}
Solving this differences inequality, we find the required uniform $L^\infty$ estimate for 
proving the convergence of time interpolators. Indeed, observe that this uniform estimate
also implies a uniform estimate for $\nor{\grad u_i^{n+1}}_2$, since by assumption $g(\bu^{n+1})$ 
remains bounded away from zero for all $n\in\N$.

\smallskip

\no\emph{Time interpolators and passing to the limit $\tau\to0$.} This step is somehow standard, so 
we give an sketch. We define, for $(t,x)\in (t_n,t_{n+1}]\times \O$, the piecewise constant and
piecewise linear interpolators
\begin{align*}
 u_i^{(\tau)}(t,x)=u_i^{n+1}(x),\quad  
 \tilde u_i^{(\tau)}(t,x)=u_i^{n}(x)+\frac{t-t_{n}}{\tau}(u_i^{n+1}(x)-u_i^{n}(x)).
\end{align*}
Using the uniform estimates of $\nor{u_i^{n+1}}_\infty$ and $\nor{\nabla u_i^{n+1}}_2$, we deduce the corresponding uniform estimates for $\nor{ u_i^{(\tau)}}_{L^{\infty}(Q_{T})}$,
$\nor{ \tilde u_i^{(\tau)}}_{L^{\infty}(Q_{T})}$, and for
\begin{equation*}
\nor{ \nabla u_i^{(\tau)}}_{L^{2}(Q_{T})},\quad  
\nor{ \nabla \tilde u_i^{(\tau)}}_{L^{2}(Q_{T})} \quad
\nor{ \p_t \tilde u_i^{(\tau)}}_{L^{2}(0,T; (H^1(\O))')}, 
\end{equation*}
 implying the existence of $ u_i\in L^2(0,T; H^1(\O))\cap L^\infty(\O)$
 and
 \[
\tilde u_i\in L^2(0,T; H^1(\O))\cap H^1(0,T; (H^1(\O))')\cap L^\infty(\O)
 \]
such that, at least in a subsequence (not relabeled), as $\tau \to 0$,
\begin{align}
 & u_i^{(\tau)} \to u_i \qtext{weakly in }L^2(0,T; H^1(\O)), \label{conv.1} \\
 & u_i^{(\tau)} \to u_i \qtext{weakly* in }L^\infty(Q_T)), \nonumber \\
 & \tilde u_i^{(\tau)} \to \tilde u_i \qtext{weakly in }   H^{1}(0,T; (H^1(\O))'), \label{conv.22}\\
 & \tilde u_i^{(\tau)} \to \tilde u_i \qtext{weakly in }   L^2(0,T; H^1(\O)), \label{conv.23}\\
 & \tilde u_i^{(\tau)} \to \tilde u_i \qtext{weakly* in }L^\infty(Q_T)).  \nonumber
\end{align}
In particular, by compactness 
\begin{equation*}
 \tilde u_i^{(\tau)} \to \tilde u_i \qtext{strongly in } L^2(Q_T). 
\end{equation*}
Since
$\tilde u_i^{(\tau)} - u_i^{(\tau)} = \tau(\delta t -1) \p_t \tilde u_i^{(\tau)}$, with $\delta t = t/\tau -n \in [0,1]$, we also find
\begin{equation}
\label{conv.25}
 \nor{\tilde u_i^{(\tau)} - u_i^{(\tau)}}_{L^2(0,T; (H^1(\O))')} \to 0\qtext{as }\tau\to 0.
\end{equation}
Hence, $u_i = \tilde u_i$. Finally, by interpolation,
we deduce 
\begin{align*}
&\nor{u_i^{(\tau)}-u_i}_{L^1(0,T; L^2(\O))} && \leq 
\nor{u_i^{(\tau)}-\tilde u_i^{(\tau)}}_{L^1(0,T; L^2(\O))}
+\nor{\tilde u_i^{(\tau)}-u_i}_{L^1(0,T; L^2(\O))}& \\
&&& \leq
\nor{ u_i^{(\tau)} - \tilde u_i^{(\tau)}}_{L^1(0,T; (H^1(\O))')}^{1/2}
\nor{ u_i^{(\tau)} - \tilde u_i^{(\tau)}}_{L^1(0,T; H^1(\O))}^{1/2} &\\
&&& +\nor{\tilde u_i^{(\tau)}-u_i}_{L^1(0,T; L^2(\O))}  \to 0,&
\end{align*}
as $\tau\to 0$, by \fer{conv.23}-\fer{conv.25}. Thus, we deduce
\begin{equation}
 \label{conv.2}
 u_i^{(\tau)} \to  u_i \qtext{strongly in } L^2(Q_T). 
\end{equation}
Considering the shift operator $\sigma_\tau u_i^{(\tau)}(t,\cdot) = u_i^{n}$, 
we  rewrite the weak form of \fer{eq:pde:ss} as 
\begin{align*}
 \int_0^T  < \p_t \tilde u_i^{(\tau)}, \vfi> + \int_{Q_T} J_i(\bu^{(\tau)})\cdot \grad \vfi =  \int_{Q_T} f_i(\bu^{(\tau)}) , 
\end{align*}
for all $\vfi \in L^2(0,T;H^1(\O))$, with $<\cdot,\cdot>$ denoting the duality product of $(H^1(\O))'\times H^1(\O)$. Finally, we pass to the limit using \fer{conv.1}, \fer{conv.22}, and \fer{conv.2}, obtaining a weak solution 
of \fer{eq:pde}-\fer{eq:bc} (with $a_{12}=a_{21}=0$).

\smallskip

\no\emph{Case 2: $a_{12}a_{21}>0$, or
 $a_{12}a_{21}=0$ and $a_{11}\neq a_{22}$. }
This case may be reduced to the diagonal case through a change of 
unknown. The eigenvalues of the matrix $A$ are given by 
 \[
  \mu_{1,2} = \frac{1}{2}\Big(\tr(A) \pm \sqrt{\Delta}\Big),\qtext{with } \Delta =(a_{11}-a_{22})^2+4a_{12}a_{21}.
 \]
Thus, the conditions on the coefficients imply $\mu_1\neq \mu_2$ (and positive). Therefore, $A$ is diagonalizable through
 a change of basis $P$, that is, there exist non-singular matrices $D$ and $P$, with $D$ diagonal, such that
 $A=P^{-1}DP$. Introducing the new unknown $\bv = P\bu$, we may write problem 
\fer{eq:pde}-\fer{def:reaction} as 
 \begin{align}
& \pt v_i- \mu_i \Div \big( g(P^{-1}\bv) \grad v_i \big) = f_i(\bv) && \qtext{in }Q_T=(0,T)\times\O, 
  \label{eq:pde2}\\
& g(P^{-1}\bv) \grad v_i\cdot \nu =0 && \qtext{on }\Gamma_T=\partial(0,T)\times\O,
\label{eq:bc2}\\
& v_{i}(\cdot,0)=v_{i0} && \qtext{in }\O, \label{eq:id2} 
\end{align}
with $\bv_0 = P \bu_0$, $ f(\bv) = B (\bv_0 - \bv)$, and $B=\diag(\{\beta_1,\beta_2\})$.

Therefore, the change of unknown renders the problem to a diagonal problem similar
to the studied in Case 1. We thus get a solution $\bv \in L^2(0,T;H^1(\O))\cap H^1(0,T;(H^1(\O))')\cap L^\infty(Q_T)$ of \fer{eq:pde2}-\fer{eq:id2}. Showing that $\bu = P^{-1}\bv$ satisfies the original problem
\fer{eq:pde}-\fer{eq:id} is immediate. $\Box$

\section{Conclusions}
We investigated the existence of solutions of a cross-diffusion model which arises as
a generalization of the image denoising model introduced by Gilboa et al. \cite{Gilboa2001}. 

For the time-independent problem with fidelity terms, we showed the existence of $L^\infty$ bounded solutions under rather general conditions on the data problem. 
This result allowed to prove the well-posedness of a quasi-steady state approximation of the evolution problem with the same  generality on the coefficient conditions as for the time-independent problem. This well-posed time-discretization of the evolution problem 
was the base for our numerical experiments. 

However, the $L^\infty$ bounds necessary to ensure the ellipticity of the diffusion operator are not directly translated from the time-independent problem to the corresponding evolution problem. Only under some rather restrictive conditions on the coefficients the passing to the limit from the QSS approximation to the evolution problem could be achieved.

Regarding the denoising capability of the method, the general conclusion is that its performance on natural images is similar to that of the gradient based Perona-Malik equation and to the Bilateral filter, and superior than the Laplacian based Perona-Malik equation. In these experiments, the Nonlocal Means filter always gave the best results. 

For textured images, the cross-diffusion method is  comparable or superior to the other methods. However, the quality gaining for these kind of images is rather poor and, probably, other more specific methods should be employed in this case. 

In our experiments with the cross-diffusion model we have always fixed the diffusion matrix as suggested in \cite{Gilboa2001, Gilboa2004}. Future work will be devoted to investigate whether other elections of the diffusion matrix coefficients   may lead to improvements of the results herein reported.

\end{document}